\documentclass[default,iicol]{sn-jnl}

\usepackage{graphicx}%
\usepackage{multirow}%
\usepackage{amsmath,amssymb,amsfonts}%
\usepackage{amsthm}%
\usepackage{mathrsfs}%
\usepackage[title]{appendix}%
\usepackage{xcolor}%
\usepackage{textcomp}%
\usepackage{manyfoot}%
\usepackage{booktabs}%
\usepackage{algorithm}%
\usepackage{algorithmicx}%
\usepackage{algpseudocode}%
\usepackage{listings}%
\usepackage{diagbox}

\newcommand{\cl}{\textcolor{red}}

\theoremstyle{thmstyleone}%
\newtheorem{theorem}{Theorem}
\theoremstyle{thmstyletwo}%
\newtheorem{lemma}{Lemma}
\theoremstyle{thmstylethree}%
\newtheorem{definition}{Definition}%

\begin{document}

\title[Article Title]{An adaptive weighted self-representation method for incomplete multi-view clustering}


\author[1]{\sur{Lishan Feng}}\email{lsfeng2022@163.com}

\author[2]{\sur{Guoxu Zhou}}\email{gx.zhou@gdut.edu.cn}

\author*[3]{\sur{Jingya Chang}}\email{jychang@gdut.edu.cn}

\affil[1]{\orgdiv{School of Mathematics and Statistics}, \orgname{Guangdong University of Technology}, 
\city{Guangzhou}, \postcode{510520}, 
\country{China}}

\affil[2]{\orgdiv{School of Automation}, \orgname{Guangdong University of Technology},
\city{Guangzhou}, \postcode{510006}, 
\country{China}}

\affil*[3]{\orgdiv{School of Mathematics and Statistics}, \orgname{Center for Mathematics and Interdisciplinary Sciences, Guangdong University of Technology}, 
\city{Guangzhou}, \postcode{510520}, 
\country{China}}


\abstract{
For multi-view data in reality, part of its elements may be missing because of human or machine error. Incomplete multi-view clustering (IMC) clusters the incomplete multi-view data according to the characters of various views of the instances. Recently, IMC has attracted much attention and many related methods have been proposed. However, the existing approaches still need to be developed and innovated in the following aspects: (1) Current methods only consider the differences of different views, while the different influences of instances, as well as distinguishes between missing values and completed values are ignored. (2) The updating scheme for weighting matrix in adaptive weighted algorithms usually relies on an optimization sub-problem, whose optimal solution may not be easy to achieve. (3) The adaptive weighted subspace algorithms that can recover the incomplete data are anchor types. The randomness of the anchor matrix may cause unreliability. To tackle these limitations, we propose an adaptive weighted self-representation (AWSR) subspace method for IMC. The AWSR method tunes the weighting matrix adaptively in accordance with the views of different instances and the recovery process of the missing values. The low rank and smoothness constraints on the representation matrix make the subspace reveal the underlying features of the dataset accurately. We also analyze the convergence property of the block coordinate method for our optimization model theoretically. Numerical performance on five real-world data shows that the AWSR method is effective and delivers superior results when compared to other eight widely-used approaches considering the clustering accuracy (ACC), normalized mutual information (NMI) and Purity.}




\maketitle

\section{Introduction}\label{sec1}


The rapid development of science and technology enables us many ways to get the information of an object. One can receive news from different news organizations or media. The same semantic can be expressed by different languages via a smart phone application program. However, in reality, due to human or machine error, such as lost files, anonymous purposes, equipment malfunctions, technical failures and so on,  it is common that the collected information is missing. Data missing can be divided into two categories: value missing, and view missing. Value missing means some of the elements of the features are absent \cite{chao2022incomplete}. View missing indicates the complete absence of some features  \cite{liu2021self, chao2022incomplete}, which is a special case of value missing. 
In this paper, we concentrate on view missing incomplete multi-view clustering (IMC) problem.

\subsection{Methods from MVC (multi-view clustering) to IMC and the auto-weighted techniques}
In this subsection, we review the kinds of subspace methods for MVC and the approaches for dealing with missing data in IMC. At last, we  summarize the auto-weighted techniques used for MVC.

Among methods of MVC, the class of subspace clustering methods is widely studied and used \cite{li2022high,zhang2017latent,fang2023comprehensive,zhuge2017robust}. These methods first map the original data into a low dimensional subspace, and then obtain the clustering results by applying spectral clustering approach to the subspace. One kind of subspace clustering approaches  assumes that the features of items in the dataset are relevant and each item can be represented by other items numerically. In this context, the self-representation subspace clustering methods arise \cite{liu2021self,tang2022consistent}. Another kind of subspace clustering methods for MVC is the anchor based approach, which chooses a small set of points called anchors as the basis of the feature space and generates the subspace of the given data by the anchors \cite{he2023scalable,sun2021scalable,guo2023scalable}. Moreover, Zhao et al. \cite{zhao2022multi} proposed the Reinforced Tensor Graph Neural Network (RTGNN) method to learn multi-view graph data, which employed tensor decomposition to obtain the graph structure features of each view in the common feature spaces. 

Compared to MVC, the missing data in IMC creates difficulties for clustering. One way of dealing with the missing data is ignoring it or simply filling it with a constant. Deng et al. introduced a missing index matrix to utilize the completed instances and  formulated a projection model for IMC \cite{deng2023projective}. Hu et al. removed the effect of missing data by assigning $0$ to it \cite{hu2021complete}. Zhao et al. suggested to use the average vector of non-missing views to estimate the missing instances \cite{zhao2022incomplete}. This kind of methods fails to explore the implicit information in the missing data and may lead to an inaccurate clustering result. The other way of dealing with the missing data is to recover it. A popular technique of restoring the missing data is to mix the matrix completion and clustering task together throughout the clustering process. Yin et al. reconstructed the incomplete data by adding a regularization term to keep the estimated matrix close enough to the existing data \cite{yin2023incomplete}. Liu et al. distinguished the complete and incomplete views and obtained the missing data via the self-representation model \cite{liu2021self}.

Different views may have different influences on the clustering result. For example, the view collected from the color of a leaf is important to the view given by the shape of a leaf when clustering data of leaves. Equal treatment of various features may bring negative effect to clustering. Therefore, weighting matrix is assigned to adjust the importance of different views \cite{shi2022self,fang2023comprehensive}. Zhao et al. proposed to use a constant diagonal matrix to measure the importance of different views \cite{zhao2022incomplete}. A majority of the weighted MVC methods are auto-weighted. The diagonal elements of the weighting matrix are restricted to the unit interval with summation being one. For subspace MVC methods, the weight is commonly multiplied by the representation error of the corresponding view \cite{wen2019unified,9914468,shu2022self,sun2021scalable,he2023scalable,liu2023auto,wan2023auto,liang2022incomplete,liu2023adaptive,khan2023adaptive}. Then the weight is updated by solving a constrained sub optimization problem whose stationary point always has a closed form. In another case, a Frobenius norm regularization term of weights is added to smooth the weight distribution \cite{ren2018robust,zhao2023auto}. In this case, the sub optimization problem  for updating the weights  becomes complex.

\subsection{Motivation and contribution}
Although great progress has been made in the area of IMC, current auto-weighted subspace IMC methods still face challenges. 1) The mechanism of weighting multiplies the weight by total representation error of the corresponding view. The features of different instances are treated equally. However, the significance of different features varies from instance to instance. For example, the color feature plays different roles in clustering leaves and persons. Also, both the known and unknown elements in the same feature are not distinguished and treated equally. 2) The adaptive updating scheme of the weighting matrix relies on the optimal solution of a constrained optimization problem, the global solution of which is not easy to find.
3) The existing auto-weighted subspace IMC methods that can recover the incomplete data are anchor based methods. The randomness of the anchor matrix induces unreliability of the clustering results.

\begin{figure*}[htpb]
	\centering
	\includegraphics[width=0.97\textwidth]{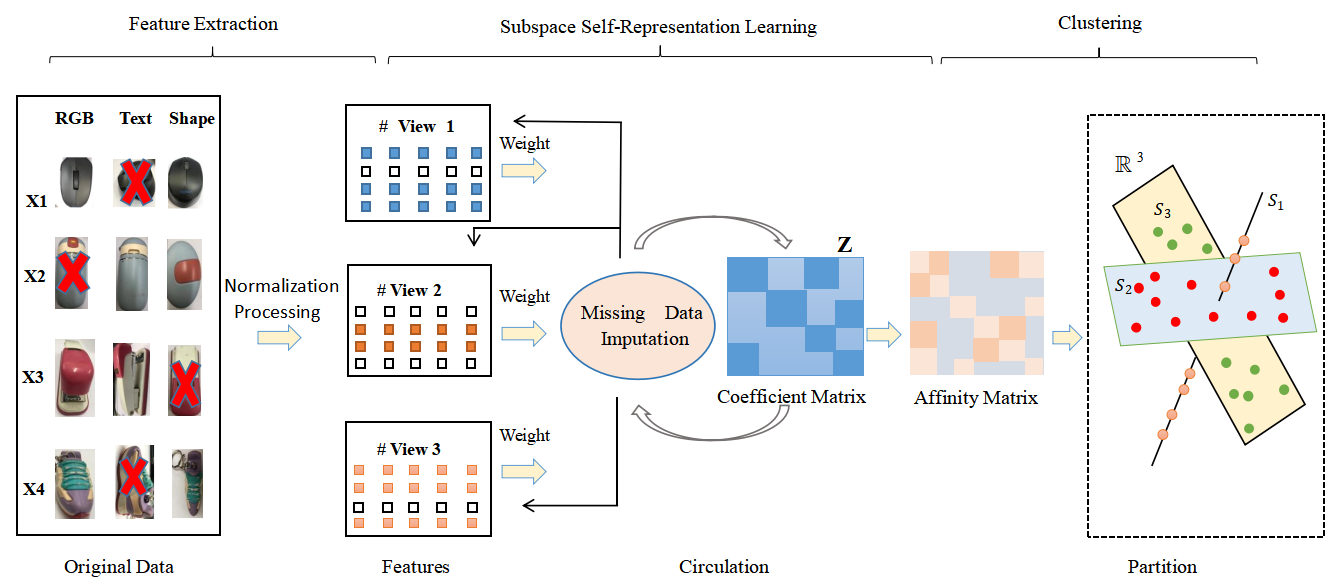}
	\caption{Overview of the algorithm. Here is an example of three views, including RGB, Text and Shape modalities. From left to right, first, feature vectors are extracted from the incomplete multi-view data $\mathbf{X}$. Then the low rank and smoothness properties are imposed on the consensus coefficient matrix $\mathbf{Z}$. Meanwhile, auto-weighted self-representation learning and missing data imputation are jointly performed to obtain the consensus representation matrix. Finally, cluster assignments are given based on the resulting affinity matrix.\label{fig1}}
\end{figure*}

In order to address the above limitations, we propose an adaptive weighted self-representation subspace clustering (AWSR) method for the IMC problem. The framework of this approach is intuitively shown in Figure \ref{fig1}.
In the AWSR model, the construction of the consensus representation matrix is designed in conjunction with the process of the missing data imputation. The weighting matrix is adaptively updated in accordance with the values assigned to the missing data and the features of different instances. Unlike traditional methods, we also put low rank and smoothness constraints on the representation matrix. The representation matrix, missing data and weighting matrix are mutually updated with each other to achieve the  optimum point.

In terms of the IMC problem, the contributions of our paper include:

(1) To the best of our knowledge, the proposed AWSR method is one of the first attempts to employ self-representation subspace algorithm with auto-weighted matrix that can obtain the representation matrix and recover the incomplete data simultaneously. A novel updating scheme of the weighting matrix is given, which not only avoids solving the sub optimization problem but also treats the values of missing data and existing data differently. Both the low rank and smoothness properties of the representation matrix are considered.

(2) In the computation process, we introduce a splitting variable which enhances the convexity of the objective function in the optimization model and makes the model easy to be solved by the block coordinate descent (BCD) algorithm. In theory, we prove the global convergence property of the iteration sequence.

(3) The numerical experiments show  complete superiority of our method over other eight state of the art approaches for IMC on five benchmark datasets.

\section{Notations and related work}

\begin{table}[htbp]
	\caption {Summary of the Notations}
	\begin{tabular}{@{}ll@{}}
		\toprule
        Notation & Description \\
        \midrule
		$\mathbf{X}_i$& The data matrix for the $i$th view.\\
		$\mathbf{Z}$ & The consensus coefficient matrix for all the views.\\
		$\mathbf{W}^{(i)}$ & The weighting matrix for the $i$th view.\\
		$\mathbf{Z}_i$ & The coefficient matrix for the $i$th view.\\
		$v$ & The number of the views.\\
		$\|\cdot\|_F$ & The Frobenius norm of matrix.\\
		$\|\cdot\|_*$ & The nuclear norm of matrix.\\
        $Tr(\cdot)$ & The trace of matrix.\\
        $\mathbf{I}$ & The identity matrix.\\
        $\varOmega$ & The feasible regions of the variables.\\
        $\otimes$ & The Kronecker product.\\
        \text{diag($\mathbf{Z}$)=0} &The diagonal elements of matrix $\mathbf{Z}$ are zero.\\
		\botrule
\end{tabular}
\end{table}
First, we introduce the notations and the preliminary knowledge. In this paper, matrix is expressed by bold capital letters, while bold lower case letters and lower case letters represent vectors and scalars respectively. The nuclear norm of a matrix $\mathbf{Z}$ is defined as $\left\|\mathbf{Z}\right\|_* = \sum\limits_i^n \sigma_i$ and $\sigma_i$ is the $i$th singular value of $\mathbf{Z}$. The Frobenius norm of $\mathbf{Z}$ is $\left\|\mathbf{Z}\right\|_F = (\sum\limits_{i=1}^n\sum\limits_{j=1}^n\mathbf{Z}_{ij}^2)^{\frac{1}{2}}.$

\begin{definition}[Subspace clustering \cite{10.1007/978-3-642-33786-4_26}]
Considering a collection of samples $\mathbf{X}=\left[\mathbf{x}_{1}, \cdots, \mathbf{x}_{n}\right].$ Suppose that $\mathbf{X}$ is derived from a combination of $k$ subspaces $\{S_i\}_{i=1}^k$ with  dimensionality $\{D_i\}_{i=1}^k$. The objective of subspace clustering is to segment samples depending on the underlying subspaces under which they are drawn from.
\end{definition}
The subspace here can be described as $S_i=\{x\in \mathbb{R}^{D_i}: x= \nu_i+ B_i y \},$ where $\nu_i$ is an arbitrary point in $S_i,$ $B_i$ is the basis of $S_i,$  and $y$ is a low-dimensional representation for $x$ \cite{vidal2011subspace}.  For example, the data points in the partition part of Figure \ref{fig1} belong to three subspaces of the $\mathbb{R}^3$ space.


One of the most important subspace clustering approach is self-representation subspace method \cite{liu2021self}. When the number of instances is large, the data of one instance can be represented by the data of others. The representation matrix represents the  latent subspace. As a surrogate of the original dataset, the latent subspace usually inherits some properties from the original data and certain constraints or properties, such as sparsity, low rank or smoothness are usually imposed on the representation matrix.
The task of self-representation subspace multi-view clustering  is  to solve
\begin{equation}
\begin{aligned}
\min\limits_{\mathbf{Z}_*,\{\mathbf{Z}_i\}_{i=1}^v} \ & \sum\limits _{i=1}^v \ell(\mathbf{X}_i-\mathbf{X}_i \mathbf{Z}_i)+\lambda \mathbf{\mathcal{R}}(\mathbf{Z}_*,\{\mathbf{Z}_i\}_{i=1}^v), \\&\text{s.t.}\ \varOmega(\mathbf{Z}_*,\{\mathbf{Z}_i\}_{i=1}^v),
\end{aligned}
\end{equation}
where $\ell$ is some loss function to measure the representation error. $\mathbf{\mathcal{R}}$ is the regularization term and  $\mathbf{Z}_*$ is the consensus coefficient matrix of all views. $\varOmega$ defines feasible regions of the variables. The regularization term is usually determined by the potential properties of the subspace.

For complete single-view dataset, if the subspace is low rank, then
the model of low rank self-representation is \cite{liu2012robust}
\begin{equation}
\min  \  \left\|\mathbf{Z}_1\right\|_*,\quad  \text{s.t.}\ \mathbf{X}_1=\mathbf{X}_1\mathbf{Z}_1, \ \text{diag}(\mathbf{Z}_1) =0.
\end{equation}
In order to capture the within-cluster affinities and smoothness properties better, the self-representation model becomes \cite{10.1007/978-3-642-33786-4_26}
\begin{equation}\label{zf2}
\begin{aligned}
\min\limits_{\mathbf{Z}_1} \ \left\|\mathbf{X}_1-\mathbf{X}_1\mathbf{Z}_1\right\|_F^2+\lambda \left\|\mathbf{Z}_1\right\|_F^2, \quad \text{s.t.}\ \text{diag}(\mathbf{Z}_1)=0.
\end{aligned}
\end{equation}


Following the framework of (Eq.\eqref{zf2}), the complete multi-view subspace clustering problem can be converted into
\begin{equation}\label{mzf3}
\begin{aligned}
&\min\limits_{\mathbf{Z}_i} \ \frac{1}{v}\sum\limits_{i=1}^v\left\|\mathbf{X}_i-\mathbf{X}_i \mathbf{Z}_i\right\|_F^2+\lambda \left\|\mathbf{Z}_i\right\|_F^2, \\
& \text{s.t.}\ \text{diag}(\mathbf{Z}_i)=0.
\end{aligned}
\end{equation}



\section{The proposed AWSR model}

When data information is incomplete, the incomplete data matrix is explicitly expressed and involved in the optimization process. Liu et al. proposed the following model \cite{liu2021self}:
\begin{equation}
\begin{aligned}
&\min\limits_{\mathbf{Z},{\{\mathbf{X}_i^{(m)}\}}_{i=1}^v,\mathbf{F}} \ \frac{1}{v}\sum\limits_{i=1}^v \left\|[\mathbf{X}_i^{(o)},\mathbf{X}_i^{(m)}]-[\mathbf{X}_i^{(o)},\mathbf{X}_i^{(m)}]\mathbf{Z}\right\|_F^2\\
&\hspace{2cm}+\lambda \left\|\mathbf{Z}\right\|_F^2 +\gamma\left\|\mathbf{Z}-\mathbf{F}^\top \mathbf{F}\right\|_F^2,\\
&\hspace{1cm} \text{s.t.}\qquad \text{diag}(\mathbf{Z})=0,\ \mathbf{F}\mathbf{ F}^\top = \mathbf{I}_k,
\end{aligned}
\end{equation}
where $\mathbf{X}_i^{(m)}$ is the  matrix saving the missing data of the $i$th view and $\mathbf{X}_i^{(o)}$ is the matrix for existing data of the $i$th view. In this model, all views share a consensus coefficient matrix $\mathbf{Z}$.


Since different views exert influences on the clustering result differently, it is a trend to include self- or auto-weighted techniques in MVC models. Two general auto-weighted  models are widely used \cite{he2023scalable,9914468,ren2018robust,zhao2023auto}:
\begin{equation}
\begin{aligned}
&\min\limits_{\mathbf{Z}_*,\{\mathbf{Z}_i\}_{i=1}^v} \  \sum\limits _{i=1}^v \text{w}_i\ell(\mathbf{X}_i-\mathbf{X}_i \mathbf{Z}_i)+\lambda \mathbf{\mathcal{R}}(\mathbf{Z}_*,\{\mathbf{Z}_i\}_{i=1}^v), \\
&\qquad \text{s.t.} \qquad   \varOmega(\mathbf{Z}_*,\{\mathbf{Z}_i\}_{i=1}^v,\textbf{W}),
\end{aligned}
\end{equation}
and
\begin{equation}
\begin{aligned}
&\min\limits_{\mathbf{Z}_*,\{\mathbf{Z}_i\}_{i=1}^v} \  \sum\limits _{i=1}^v \text{w}_i\ell(\mathbf{X}_i-\mathbf{X}_i \mathbf{Z}_i)+\lambda_1 \mathbf{\mathcal{R}}(\mathbf{Z}_*,\{\mathbf{Z}_i\}_{i=1}^v), \\
&\hspace{1.5cm} + \lambda_2\mathbf{\mathcal{R}}(\textbf{W}) \\
&\qquad \text{s.t.} \qquad   \varOmega(\mathbf{Z}_*,\{\mathbf{Z}_i\}_{i=1}^v,\textbf{W}).
\end{aligned}
\end{equation}
The vector $\textbf{w}=(\text{w}_1,\ldots,\text{w}_v)$ is composed of the diagonal entries of the diagonal weighting matrix $\textbf{W}$. The regularization term $\mathbf{\mathcal{R}}(\textbf{W})$ is the Frobenius norm in \cite{ren2018robust,zhao2023auto}.

As pointed in \cite{fang2023comprehensive}, when coping with various real-life applications, auto-weighting should be more adaptive and automatic. The existing weighting matrix only concentrate on the differences among influences of different views. However the differences of instances, as well as dissimilarities between the missing values and existing values are neglected. Besides, the weighting matrix is updated by solving a sub optimization problem, which may be troublesome. To overcome these difficulties, we suggest a new way to automatically produce the weighting matrix. Initially the weighting matrix is a diagonal matrix with its diagonal elements being
 \begin{equation}
\mathbf{W}_{jj}^{(i)}=\left\{
\begin{array}{cl}
1 , & \mbox{if $j$th item is in the $i$th view,}\\
0 , & \text{otherwise}.\\
\end{array}\right.
\end{equation}
Because the missing view is totally unknown at first, the zero diagonal elements in this weighting matrix make the representation error caused by the corresponding missing view being zero. At the beginning, the value of the data in the missing view is set to zero. During the iteration process,  some elements in the missing view of the item are recovered while other elements are still zero. In this case, the influence of the recovered data should be taken into account. Therefore, the elements of the weighting matrix at the index positions where the missing data is recovered are set to
\begin{equation}
\mathbf{W}_{mj}^{(i)}=\left\{
\begin{array}{cl}
\alpha_{mj}^{(i)} , & \mbox{if $\mathbf{X}_{i,mj}$ is recovered or nonzero}\\
0 , & \text{otherwise}.
\end{array}\right.
\end{equation}
The updating scheme means if the $m$th element of the $j$th instance in the $i$th missing view $\mathbf{X}_{i,mj}$ is recovered, then $\mathbf{W}_{mj}^{(i)}=\alpha^i_{mj}.$ Here $\alpha^i_{mj}$ is nonzero and determined by either the value of the missing element $\mathbf{X}_{i,mj}$ or an empirical constant. On the other hand, if the $n$th element of the $j$th instance in the $i$th missing view $\mathbf{X}_{i,nj}$ is zero, then $\mathbf{W}_{nj}^{(i)}=0.$ We illustrate the influence tuned by $\mathbf{W}_{mj}^{(i)}$ and $\mathbf{W}_{nj}^{(i)}$ in detail. Let $\mathbf{E}^i = \mathbf{X}_i-\mathbf{X}_i \mathbf{Z}.$ The matrix product
 $$
\mathbf{E}^i (\mathbf{W}^{(i)})^{\top}=\left(
  \begin{array}{ccccc}
    \cdots & \mathbf{E}^{i}_{1j} & \cdots \\
     & \vdots &  \\
    \cdots & \mathbf{E}^{i}_{pj} & \cdots \\
  \end{array}
\right)
\left(
  \begin{array}{ccccc}
    \cdots & \vdots & \cdots & \vdots & \cdots \\
    \cdots & \mathbf{W}_{mj}^{(i)} &\cdots &  \mathbf{W}_{nj}^{(i)}&  \cdots \\
    \cdots &\vdots & \cdots & \vdots & \cdots\\
  \end{array}
\right)
$$
implies that $\mathbf{W}_{mj}^{(i)}=\alpha^i_{mj}$ affirms the effect of the $j$th instance in the $i$th view on clustering results, while $\mathbf{W}_{nj}^{(i)}=0$ weakens the effect of the $j$th instance in the $i$th view on clustering results.
Therefore, by utilizing the adaptive weighting matrix technique we achieve the goal that the recovered data has significance on the clustering result and the significance of the unrecovered data on clustering result is eliminated. The influences of different instances in the same view are distinguished. Furthermore, the updating framework avoids solving the optimization sub-problem of $\mathbf{W}^{(i)}$.

In order to get a low rank and smoothness coefficient matrix, the nuclear norm and Frobenius norm are employed to construct regularization terms, which can also capture the data correlation structure in the same space and help to solve the minimization problem effectively. Finally, the overall model is
\begin{equation} \label{start obj}
\begin{aligned}
\min\limits_{\mathbf{Z},{\{\mathbf{X}_i^{(m)}\}}_{i=1}^v} \ &\sum\limits_{i=1}^v\left\|([\mathbf{X}_i^{(o)},\mathbf{X}_i^{(m)}]-[\mathbf{X}_i^{(o)},\mathbf{X}_i^{(m)}]\mathbf{Z}){\mathbf{W}^{(i)}}^\top\right\|_F^2\\&+\gamma\left\|\mathbf{Z}\right\|_*+ \frac{\lambda}{2} \left\|\mathbf{Z}\right\|_F^2 ,\\
&\text{s.t.} \ \text{diag}(\mathbf{Z})=0,
\end{aligned}
\end{equation}
where
$\gamma, \lambda$ are positive trade-off parameters and $\mathbf{X}_i^{(o)}$ and $\mathbf{X}_i^{(m)}$ refer to observed and missing entries of the $i$th data view. 

%
%
%
%
%
%


In summary, the entire framework of our method are as follows. First we construct the adaptive weighted linear constrained optimization model (Eq.\eqref{start obj}) based on the subspace representation approach. Then, problem (Eq.(\ref{start obj})) is solved  iteratively by the BCD method. Finally, the Laplacian matrix is generated based on the coefficient matrix $\mathbf{Z}$ and the spectral clustering approach \cite{CAI2018316} is employed to  cluster the given objects.

\section{Computation of the AWSR model}
 In this section, we demonstrate computation process of the optimization model (Eq.\eqref{start obj}). At first, we declare that the nonzero elements of the weighting matrix W is determined artificially according to the numerical results.
In order to enhance the convex property of the objective function and create a separable optimization model, we introduce a separated variable $\mathbf{J}$ and convert model (Eq.(\ref{start obj})) to
\begin{equation}\label{all obj function}
\begin{aligned}
&\min\limits_{\mathbf{J},\mathbf{Z},{\{\mathbf{X}_i^{(m)}\}}_{i=1}^v}\ \gamma\left\|\mathbf{Z}\right\|_*+\frac{\lambda}{2} \left\|\mathbf{Z}\right\|_F^2 + \frac{\alpha}{2} \left\|\mathbf{J}-\mathbf{Z}\right\|_F^2\\&+\sum\limits_{i=1}^v\left\|([\mathbf{X}_i^{(o)},\mathbf{X}_i^{(m)}]-[\mathbf{X}_i^{(o)},\mathbf{X}_i^{(m)}]\mathbf{J}){\mathbf{W}^{(i)}}^\top\right\|_F^2,\\
& \text{s.t.}  \  \text{diag}(\mathbf{Z})=0.
\end{aligned}
\end{equation}
The (Eq.(\ref{start obj})) and (Eq.(\ref{all obj function})) are equivalent when $\alpha > 0$ is large enough. The strong convexity property of the objective function implies that the solution is stable and unique.
We apply the alternate strategy \cite{liu2021hierarchical} to calculation of the above optimization problem.

\subsection{Sub-problem for J:} To update $\mathbf{J}$, we fix the values of all other variables. The sub-problem related with $\mathbf{J}$ is
\begin{equation} \label{eq15}
\begin{aligned}
\min\limits_{\mathbf{J}}  \psi(\mathbf{J})  &= \sum\limits_{i=1}^v\left\|([\mathbf{X}_i^{(o)},\mathbf{X}_i^{(m)}]-[\mathbf{X}_i^{(o)},\mathbf{X}_i^{(m)}]\mathbf{J}){\mathbf{W}^{(i)}}^\top\right\|_F^2\\&+\frac{\alpha}{2} \left\|\mathbf{J}-\mathbf{Z}\right\|_F^2.
\end{aligned}
\end{equation}
The (Eq.(\ref{eq15})) is equivalent to the following problem
\begin{equation} \label{4.1-4}
\begin{aligned}
\min\limits_{\mathbf{J}} \ \psi(\mathbf{J}) & = Tr( \sum\limits_{i=1}^v \mathbf{X}_i \mathbf{J} {\mathbf{W}^{(i)}}^\top \mathbf{W}^{(i)} \mathbf{J}^\top {\mathbf{X}_i}^\top \\&-\mathbf{X}_i \mathbf{J} {\mathbf{W}^{(i)}}^\top \mathbf{W}^{(i)} {\mathbf{X}_i}^\top - \mathbf{X}_i {\mathbf{W}^{(i)}}^\top \mathbf{W}^{(i)} {\mathbf{J}}^\top \mathbf{X}_i ) \\&+ \frac{\alpha}{2} Tr(\mathbf{J}^\top \mathbf{J} -\mathbf{J}^\top \mathbf{Z}-\mathbf{Z}^\top \mathbf{J} ),
\end{aligned}
\end{equation}
in which $\mathbf{X}_i=[\mathbf{X}_i^{(o)},\mathbf{X}_i^{(m)}]$.
Take derivation of the objective function in (Eq.\eqref{4.1-4}) with respect to variable $\mathbf{J}$, and let the derivation equal to zero. We get
\begin{equation}\label{4.1-1}
\begin{aligned}
&\sum\limits_{i=1}^v {\mathbf{X}_i}^\top \mathbf{X}_i \mathbf{J} {\mathbf{W}^{(i)}}^\top \mathbf{W}^{(i)} + \frac{\alpha}{2} \mathbf{J} \\&-\sum\limits_{i=1}^v{\mathbf{X}_i}^\top \mathbf{X}_i {\mathbf{W}^{(i)}}^\top \mathbf{W}^{(i)} -\frac{\alpha}{2} \mathbf{Z} = 0.
\end{aligned}
\end{equation}
Since the matrix equation $\mathbf{AXB} = \mathbf{C}$ is equivalent to \cite{merino1992topics}
$$(\mathbf{B}^T\otimes \mathbf{A})\text{vec} (\mathbf{X}) = \text{vec}(\mathbf{C}),$$
where $\otimes$ means the Kronecker product, the (Eq.\eqref{4.1-1}) equals to
\begin{equation}\label{4.1-2}
\left[\sum\limits_{i=1}^v {\mathbf{W}^{(i)}}^\top \mathbf{W}^{(i)} \otimes {\mathbf{X}_i}^\top \mathbf{X}_i   +\frac{\alpha}{2} \mathbf{I} \right] \text{vec}(\mathbf{J}) = \text{vec}(\mathbf{C}).
\end{equation}
Here $\mathbf{C} = \sum\limits_{i=1}^v{\mathbf{X}_i}^\top \mathbf{X}_i {\mathbf{W}^{(i)}}^\top \mathbf{W}^{(i)} +\frac{\alpha}{2}\mathbf{Z}. $ The coefficient matrix in the linear equation system
(Eq.\eqref{4.1-2}) is positive definite. Hence, the vector $\text{vec}(\mathbf{J})$ can be solved by the traditional conjugate gradient method. Finally, we get the matrix $\mathbf{J}$ by rearranging $\text{vec}(\mathbf{J}).$

\subsection{Sub-problem for $\mathbf{X}_i^{(m)}$ :} To update $\mathbf{X}_i^{(m)}$, we fix the values of all other variables. The sub-problem related with $\mathbf{X}_i^{(m)}$ is
\begin{equation}\label{eq23}
\begin{aligned}	\min\limits_{\mathbf{X}_i^{(m)}} \ &\psi(\mathbf{X}_i^{(m)})=\\& \sum\limits_{i=1}^v\left\|([\mathbf{X}_i^{(o)},\mathbf{X}_i^{(m)}]-[\mathbf{X}_i^{(o)},\mathbf{X}_i^{(m)}]\mathbf{J}){\mathbf{W}^{(i)}}^\top\right\|_F^2.
\end{aligned}
\end{equation}
The (Eq.(\ref{eq23})) can be transformed to
\begin{equation}\label{eq24}
\begin{aligned}
\min\limits_{\mathbf{X}_i^{(m)}} \ &\psi(\mathbf{X}_i^{(m)})=Tr([\mathbf{X}_i^{(o)},\mathbf{X}_i^{(m)}]\mathbf{B}_i{[\mathbf{X}_i^{(o)},\mathbf{X}_i^{(m)}]}^\top) , \\& \text{for} \ i=1,\ldots,v,
\end{aligned}
\end{equation}
where
\begin{equation}
\begin{aligned}
\mathbf{B}_i &={\mathbf{W}^{(i)}}^\top \mathbf{W}^{(i)} - \mathbf{J}{\mathbf{W}^{(i)}}^\top \mathbf{W}^{(i)}- {\mathbf{W}^{(i)}}^\top \mathbf{W}^{(i)} \mathbf{J}^\top \\& + \mathbf{J}{\mathbf{W}^{(i)}}^\top \mathbf{W}^{(i)} \mathbf{J}^\top.
\end{aligned}
\end{equation}
 Because $\mathbf{X}_i$ is divided into observed and missing data, i.e., $\mathbf{X}_i=[\mathbf{X}_i^{(o)},\mathbf{X}_i^{(m)}]$, we rewrite $\mathbf{B}_i$ in the form of block matrix. The (Eq.(\ref{eq24})) is equivalent to
\begin{equation}\label{eq25}
\begin{aligned}
\min\limits_{\mathbf{X}_i^{(m)}}& \  \psi(\mathbf{X}_i^{(m)})=\\& Tr([\mathbf{X}_i^{(o)},\mathbf{X}_i^{(m)}]\left(\begin{array}{cc}
\mathbf{B}_i^{(oo)}& \mathbf{B}_i^{(om)}\\
\mathbf{B}_i^{(mo)}& \mathbf{B}_i^{(mm)}
\end{array}\right)
{[\mathbf{X}_i^{(o)},\mathbf{X}_i^{(m)}]}^\top.
\end{aligned}
\end{equation}
The (Eq.\eqref{eq25}) is further replaced by
\begin{equation} \label{4.2-1}
\begin{aligned}
&\min\limits_{\mathbf{X}_i^{(m)}}
\ \psi(\mathbf{X}_i^{(m)}) =  \\&Tr(\mathbf{X}_i^{(o)}(\mathbf{B}_i^{(om)}+{\mathbf{B}_i^{(mo)}}^\top){\mathbf{X}_i^{(m)}}^\top+\mathbf{X}_i^{(m)}\mathbf{B}_i^{(mm)}{\mathbf{X}_i^{(m)}}^\top).
\end{aligned}
\end{equation}
From the definition of $\mathbf{B}_i$, we have
\begin{equation}
\mathbf{B}_i^{(mm)} = \left(\mathbf{W}^{(i)}-\mathbf{J} \mathbf{W}^{(i)}\right) ^{(m)} \left[\left(\mathbf{W}^{(i)}-\mathbf{J} \mathbf{W}^{(i)}\right)^{(m)}\right] ^\top ,
\end{equation}
in which $(\mathbf{W}^{(i)}-\mathbf{J} \mathbf{W}^{(i)}) ^{(m)}$ refers to the columns of $\mathbf{W}^{(i)}-\mathbf{J} \mathbf{W}^{(i)}$ and corresponds to the weighted missing entries in the $i$th view.
By taking the derivation of the objective function in
(Eq.\eqref{4.2-1}) with respect to the variable $\mathbf{X}_i^{(m)}$, we get
\begin{equation}\label{4.2-2}
\mathbf{X}_i^{(m)} \mathbf{B}_i^{(mm)}=-\mathbf{X}_i^{(o)} \mathbf{B}_i^{(om)}.
\end{equation}
Because $\mathbf{B}_i^{(mm)}$ is a positive semi-definite matrix,
we calculate the eigenvalue decomposition of matrix $\mathbf{B}_i^{(mm)}$, i.e., $\mathbf{B}_i^{(mm)} = \mathbf{U}_i \mathbf{\wedge}_i \mathbf{V}_i^{\top}$.
The numerical solution of equation system (Eq.\eqref{4.2-2}) is
\begin{equation}\label{Xsolution}
\mathbf{X}_i^{(m)} =  - \mathbf{X}_i^{(o)} \mathbf{B}_i^{(om)}\mathbf{V}_i\mathbf{\wedge}_i^{+} {\mathbf{U}_i}^{\top}.
\end{equation}
Here $\mathbf{\wedge}_i^{+}$ is the pseudo inverse of $\mathbf{\wedge}_i$.

\subsection{Sub-problem for $Z$:} To update $\mathbf{Z}$, we fix the values of all other variables. We solve the following sub-problem:
\begin{equation}\label{eq28}
\begin{aligned}
&\min\limits_{\mathbf{Z}}\  \psi(\mathbf{Z})=  \  \gamma\left\|\mathbf{Z}\right\|_*+ \frac{\lambda}{2} \left\|\mathbf{Z}\right\|_F^2+\frac{\alpha}{2} \left\|\mathbf{J}-\mathbf{Z}\right\|_F^2  ,\\
&\text{s.t.} \ \ \text{diag}(\mathbf{Z})=0.
\end{aligned}
\end{equation}
The objective function in (Eq.\eqref{eq28}) is rearranged, and we obtain the equivalent model as follows
\begin{equation}\label{eq11}
\begin{aligned}
&\min\limits_{\mathbf{Z}}\ \psi(\mathbf{Z}) =\ \frac{\gamma}{\lambda+\alpha}\left\|\mathbf{Z}\right\|_* + \frac{1}{2} \left\|\mathbf{Z}-\frac{\alpha}{\lambda+\alpha} \mathbf{J}\right\|_F^2,\\
&\text{s.t.} \ \text{diag}(\mathbf{Z})=0.
\end{aligned}
\end{equation}
Denote the Lagrange function $$\mathcal{L}(\mathbf{Z},\mathbf{y}) = \psi(\mathbf{Z}) + <\mathbf{y},\text{diag}(\mathbf{Z})>.$$ Then the dual function of the Lagrange multiplier $\mathbf{y}$ is $g(\mathbf{y})=\inf\limits_\mathbf{Z} \mathcal{L}(\mathbf{Z},\mathbf{y}).$  The optimization problem (Eq.\eqref{eq11}) is convex and satisfies the weak Slater condition, so the strong duality property holds that means
\begin{equation}
\sup\limits_\mathbf{y} \inf\limits_\mathbf{Z} \mathcal{L}(\mathbf{Z},\mathbf{y}) = \mathcal{L}(\mathbf{Z}^*,\mathbf{y}^*)=\inf\limits_\mathbf{Z} \sup\limits_\mathbf{y} \mathcal{L}(\mathbf{Z},\mathbf{y}).
\end{equation}
We employ Uzawa's algorithm \cite{bramble1997analysis} to solve the dual problem of (Eq.\eqref{eq11}). The gradient of $g(\mathbf{y})$ is given by
\begin{equation}
\nabla g(\mathbf{y})=\frac{\partial \mathcal{L}(\tilde{\mathbf{Z}},\mathbf{y})}{\partial \mathbf{y}}= \text{diag}(\tilde{\mathbf{Z}}),
\end{equation}
where $\tilde{\mathbf{Z}}=\text{arg} \min\limits_\mathbf{Z}\mathcal{L}(\mathbf{Z},\mathbf{y}).$
Uzawa's algorithm is achieved iteratively via the following scheme
\begin{equation}
\begin{cases} \label{Uz}
\mathcal{L}(\mathbf{Z}^{t},\mathbf{y}^{t-1}) = \min\limits_\mathbf{Z} \mathcal{L}(\mathbf{Z},\mathbf{y}^{t-1}),\\
\mathbf{y}^{t} = \mathbf{y}^{t-1} + \delta_k (\text{diag}(\mathbf{Z}^t)).
\end{cases}
\end{equation}
The parameter $\delta_k$ is a step size of the gradient direction. Next, we consider the minimization problem $\min\limits_Z \mathcal{L}(\mathbf{Z},\mathbf{y}^{t-1})$ in (Eq.\eqref{Uz}). By substituting the formula in (Eq.\eqref{eq11}) for $\psi(\mathbf{Z}),$ we obtain
\begin{equation}\label{4.3-1}
\begin{aligned}
\mathcal{L}(\mathbf{Z},\mathbf{y}^{t-1}) &
= \psi(\mathbf{Z}) + <\mathbf{y}^{t-1},\text{diag}(\mathbf{Z})>\\
& = \psi(\mathbf{Z}) + <\mathbf{Y}^{t-1},\mathbf{Z}>.
\end{aligned}
\end{equation}
Here $\mathbf{Y}$ is a matrix of the same size as $\mathbf{Z}$. Its diagonal elements equal to $\mathbf{y}$ and off diagonal elements are zero.
Furthermore, we have the following equivalent model
\begin{equation}\label{4.3-3}
\begin{aligned}
&\text{arg}\min\limits_{\mathbf{Z}}\mathcal{L}(\mathbf{Z},\mathbf{y}^{t-1})=\frac{\gamma}{\lambda+\alpha} \left\|\mathbf{Z}\right\|_*\\
&+\text{arg}\min\limits_{\mathbf{Z}} \frac{1}{2}\left\|\mathbf{Z}-\left(\frac{\alpha}{\lambda+\alpha}(\mathbf{J}+ \frac{\mathbf{Y}^{t-1}}{\alpha})\right)\right\|_F^2.
\end{aligned}
\end{equation}
Based on (Eq.\eqref{4.3-3}), the minimization problem in (Eq.\eqref{Uz}) can be efficiently solved by the singular value threshold operator \cite{cai2010singular} and Uzawa's process is transformed into
\begin{equation}\label{Zsolution}
\begin{cases}
\begin{aligned}
\mathbf{Z}^t &=\mathcal{D}_{\frac{\gamma}{\lambda+\alpha}}\left(\frac{\alpha}{\lambda+\alpha}(\mathbf{J}+ \frac{\mathbf{Y}^{t-1}}{\alpha})\right)\\
\mathbf{y}^t &= \mathbf{y}^{t-1} + \delta_k (\text{diag}(\mathbf{Z}^t)),\\
\end{aligned}
\end{cases}
\end{equation}
where
\begin{equation*}\label{4.3-2}
\begin{aligned}
& \mathcal{D}_{\frac{\gamma}{\lambda+\alpha}}
\left(\frac{\alpha}{\lambda+\alpha}(\mathbf{J}+ \frac{\mathbf{Y}^{t-1}}{\alpha})\right)=
\mathbf{U}\mathcal{D}_{\frac{\gamma}{\lambda+\alpha}}(\mathbf{\Sigma}) \mathbf{V}^T \\
& \mathcal{D}_{\frac{\gamma}{\lambda+\alpha}}(\mathbf{\Sigma}) = \text{diag}
\{\max(0,\sigma_i -\frac{\gamma}{\lambda+\alpha})\}
\end{aligned}
\end{equation*}
and $\mathbf{U}\mathbf{\Sigma} \mathbf{V}$ is the singular value decomposition of the matrix $$\frac{\alpha}{\lambda+\alpha}(\mathbf{J}+ \frac{\mathbf{Y}^{t-1}}{\alpha}) \text{ with} \quad \mathbf{\Sigma} = \text{diag}(\sigma_i).$$

When the iteration converges, we obtain the consensus coefficient matrix $\mathbf{Z}$. Finally, we apply the spectral clustering technology \cite{CAI2018316} to the Laplacian (affinity) matrix $\mathbf{L}=(|\mathbf{Z}|+|\mathbf{Z}^\top|)/2$ to calculate the cluster assignments. The procedure for solving our model (Eq.\eqref{all obj function}) is summarized in Algorithm \ref{algorithm1}.

\begin{algorithm}
	\caption{An adaptive weighted self-representation method for incomplete multi-view clustering (AWSR)}\label{algorithm1}
	\begin{algorithmic}[1]
		
		\Require Initial weighting matrix $\mathbf{W}_0^{(i)}$, incomplete multi-view data ${\{[\mathbf{X}_i^{(o)},\mathbf{X}_i^{(m)}]\}}_{i=1}^v,$  initial matrix $\mathbf{Z}_0$, $\mathbf{J}_0,$ the number of clusters $k$, parameter $\gamma$, $\lambda$, $\alpha$, $error=1e-3$;
		
		\While{abs((obj(t-1)-obj(t))/(obj(t)))$>$ error}
		
		\State  Update $\mathbf{J}$ by solving Eq.(\ref{4.1-2});
		\State  Update $\{\mathbf{X}_i^{(m)}\}_{i=1}^{v}$ by Eq.(\ref{Xsolution});
		\State  Update $\mathbf{Z}$ by Eq.(\ref{Zsolution});
		\State  t=t+1;
		
		\EndWhile
		\State  Obtain the affinity matrix by $\mathbf{L}=(|\mathbf{Z}|+|\mathbf{Z}^\top|)/2$;
        \State  Apply spectral clustering method with the aid of the affinity matrix $\mathbf{L}$;
        \Ensure Clustering result $\mathcal{C}$.
	\end{algorithmic}
\end{algorithm}

\subsection{Convergence and complexity analysis}
In this section, we prove that if the iteration sequence generated by the  AWSR algorithm is infinite, it converges to a stationary point globally. Also, we demonstrate the convergence property of the sequences of objective function values numerically.

Let $\tilde{\mathbf{X}}^{(m)} = \{\mathbf{X}_1^{(m)}, \mathbf{X}_2^{(m)}, \cdot\cdot\cdot ,\mathbf{X}_v^{(m)}\} $. First, we demonstrate the property of $\{\mathbf{J}^{(t)},\mathbf{Z}^{(t)},\tilde{\mathbf{X}}^{(m)(t)}\}$ in the following lemma.
\begin{lemma}\label{lemma2}
	The sequence \{$\mathbf{J}^{(t)},\mathbf{Z}^{(t)},\tilde{\mathbf{X}}^{(m)(t)}$\} generated via Algorithm \ref{algorithm1} satisfies those properties:
	
	(1) The function $f(\mathbf{J}^{(t)},\mathbf{Z}^{(t)},\tilde{\mathbf{X}}^{(m)(t)})$ is a monotonically decreasing function. Especially, we have
	\begin{equation}
	\begin{aligned}
	&f(\mathbf{J}^{(t+1)},\mathbf{Z}^{(t+1)},\tilde{\mathbf{X}}^{(m)(t+1)})\\
	&\leqslant f(\mathbf{J}^{(t)},\mathbf{Z}^{(t)},\tilde{\mathbf{X}}^{(m)(t)}) - \frac{\alpha}{2}\left\|\mathbf{J}^{(t+1)}-\mathbf{J}^{(t)}\right\|_F^2\\& - \frac{\lambda}{2}\left\|\mathbf{Z}^{(t+1)}-\mathbf{Z}^{(t)}\right\|_F^2
	\end{aligned}
	\end{equation}
	
	(2) $\lim\limits_{t \rightarrow \infty}f(\mathbf{J}^{(t)},\mathbf{Z}^{(t)},\tilde{\mathbf{X}}^{(m)(t)}) = c$ for some constant $c$.
	
	(3) When $t \rightarrow +\infty$, $\mathbf{J}^{(t+1)}-\mathbf{J}^{(t)}\rightarrow 0$, $\mathbf{Z}^{(t+1)}-\mathbf{Z}^{(t)}\rightarrow 0$ and $\tilde{\mathbf{X}}^{(m)(t+1)}-\tilde{\mathbf{X}}^{(m)(t)}\rightarrow 0$.
	
	(4) The sequences $\{ \mathbf{J}^{(t)}\}$, $\{\mathbf{Z}^{(t)}\}$ and $\{\tilde{\mathbf{X}}^{(m)(t)}\}$ are bounded.
\end{lemma}

\begin{proof}
	\textbf{(1).} The function $f(\mathbf{J},\mathbf{Z}^{(t)},\tilde{\mathbf{X}}^{(m)(t)})$ of $\mathbf{J}$ in the $t$th updating process of $\mathbf{J}$ (Eq.\eqref{eq15}) is a $\alpha$-strongly convex function. Therefore, we have
	\begin{equation}\label{lemma-1}
	\begin{aligned}
	f(\mathbf{J}^{(t+1)},\mathbf{Z}^{(t)},\tilde{\mathbf{X}}^{(m)(t)})& \leqslant f(\mathbf{J}^{(t)},\mathbf{Z}^{(t)},\tilde{\mathbf{X}}^{(m)(t)})\\& - \frac{\alpha}{2}\left\|\mathbf{J}^{(t+1)}-\mathbf{J}^{(t)}\right\|_F^2.
	\end{aligned}
	\end{equation}
	The function $f(\mathbf{J}^{(t+1)},\mathbf{Z}^{(t)},\tilde{\mathbf{X}}^{(m)})$ of $\tilde{\mathbf{X}}^{(m)}$ in the $t$th updating process of $\tilde{\mathbf{X}}^{(m)}$ (Eq.\eqref{eq23}) is a convex function, and the following inequality holds
	\begin{equation}\label{lemma-2}
	\begin{aligned}
	f(\mathbf{J}^{(t+1)},\mathbf{Z}^{(t)},\tilde{\mathbf{X}}^{(m)(t+1)}) \leq f(\mathbf{J}^{(t+1)},\mathbf{Z}^{(t)},\tilde{\mathbf{X}}^{(m)(t)}).
	\end{aligned}
	\end{equation}
	Similarly, the function $f(\mathbf{J}^{(t+1)},\mathbf{Z},\tilde{\mathbf{X}}^{(m)(t+1)})$ of $\mathbf{Z}$ is a $\lambda$-strongly convex function in the sub-problem (Eq.\eqref{eq28}), and we have
	\begin{equation}\label{lemma-3}
	\begin{aligned}
	f(\mathbf{J}^{(t+1)},\mathbf{Z}^{(t+1)},\tilde{\mathbf{X}}^{(m)(t+1)})& \leqslant f(\mathbf{J}^{(t+1)},\mathbf{Z}^{(t)},\tilde{\mathbf{X}}^{(m)(t+1)}) \\&- \frac{\lambda}{2}\left\|\mathbf{Z}^{(t+1)}-\mathbf{Z}^{(t)}\right\|_F^2 .
	\end{aligned}
	\end{equation}
	Summation of the above ( Eq.\eqref{lemma-1}), (Eq.\eqref{lemma-2}) and (Eq.\eqref{lemma-3}) induces the inequality
	\begin{equation}\label{eq37}
	\begin{aligned}
	f(\mathbf{J}^{(t+1)},\mathbf{Z}^{(t+1)},\tilde{\mathbf{X}}^{(m)(t+1)})  &\leqslant f(\mathbf{J}^{(t)},\mathbf{Z}^{(t)},\tilde{\mathbf{X}}^{(m)(t)})\\&-\frac{\lambda}{2}\left\|\mathbf{Z}^{(t+1)}-\mathbf{Z}^{(t)}\right\|_F^2\\&
	- \frac{\alpha}{2}\left\|\mathbf{J}^{(t+1)}-\mathbf{J}^{(t)}\right\|_F^2,
	\end{aligned}
	\end{equation}
	which means the function $f(\mathbf{J}^{(t)},\mathbf{Z}^{(t)},\tilde{\mathbf{X}}^{(m)(t)})$ decreases as $t$ increases.
	
	\textbf{(2).} Because $f(\mathbf{J}^{(t)},\mathbf{Z}^{(t)},\tilde{\mathbf{X}}^{(m)(t)})$ is non-negative and monotonically decreasing, it converges as $t$ tends to infinity. That is to say $$\lim\limits_{t\rightarrow +\infty} f(\mathbf{J}^{(t)},\mathbf{Z}^{(t)},\tilde{\mathbf{X}}^{(m)(t)}) = c \ \text{for some constant }\ c.$$
	
	\textbf{(3).} Based on (Eq.\eqref{eq37}), we have
	\begin{equation}
	\begin{aligned}
	& \sum\limits_{t=0}^{+\infty}
	[\frac{\lambda}{2}\left\|\mathbf{Z}^{(t+1)}-\mathbf{Z}^{(t)}\right\|_F^2
	+ \frac{\alpha}{2}\left\|\mathbf{J}^{(t+1)}-\mathbf{J}^{(t)}\right\|_F^2]\\
	\leqslant& \sum\limits_{t=0}^{+\infty}[f(\mathbf{J}^{(t)},\mathbf{Z}^{(t)},\tilde{\mathbf{X}}^{(m)(t)})\\&-f(\mathbf{J}^{(t+1)},\mathbf{Z}^{(t+1)},\tilde{\mathbf{X}}^{(m)(t+1)})]\\ \leq &f(\mathbf{J}^{(0)},\mathbf{Z}^{(0)},\tilde{\mathbf{X}}^{(m)(0)}).
	\end{aligned}
	\end{equation}
	Furthermore, we get $\mathbf{J}^{(t+1)}-\mathbf{J}^{(t)}\rightarrow 0$, $\mathbf{Z}^{(t+1)}-\mathbf{Z}^{(t)}\rightarrow 0$. In addition, the rule of updating $\tilde{\mathbf{X}}^{(m)(t)}$ indicates $$f(\mathbf{J}^{(t+1)},\mathbf{Z}^{(t)},\tilde{\mathbf{X}}^{(m)(t+1)})\leq f(\mathbf{J}^{(t+1)},\mathbf{Z}^{(t)},\tilde{\mathbf{X}}^{(m)(t)}),$$ and then we obtain $\tilde{\mathbf{X}}^{(m)(t+1)}-\tilde{\mathbf{X}}^{(m)(t)}\rightarrow 0$ from (Eq.\eqref{eq24}).
	
	\textbf{(4).} Due to the fact that the value $f(\mathbf{J}^{(t)},\mathbf{Z}^{(t)},\tilde{\mathbf{X}}^{(m)(t)})$ is bounded and each term in (Eq.\eqref{all obj function}) is nonnegative, the sequences $\{ \mathbf{J}^{(t)}\}$, $\{\mathbf{Z}^{(t)}\}$ and $\{\tilde{\mathbf{X}}^{(m)(t)}\}$ are all bounded.
\end{proof}

\begin{theorem}\label{theorem1}
	Suppose $\{(\mathbf{J}^{(t)},\mathbf{Z}^{(t)},\tilde{\mathbf{X}}^{(m)(t)})\}$ is the infinite iteration sequence produced by Algorithm \ref{algorithm1}. Then there exists $(\mathbf{J}^{(*)},\mathbf{Z}^{(*)},\tilde{\mathbf{X}}^{(m)(*)})$ such that
	\begin{equation}\label{Th1}
	\lim_{t\rightarrow \infty}
	(\mathbf{J}^{(t)},\mathbf{Z}^{(t)},\tilde{\mathbf{X}}^{(m)(t)})=(\mathbf{J}^{(*)},\mathbf{Z}^{(*)},\tilde{\mathbf{X}}^{(m)(*)})
	\end{equation}
	with $\text{diag}(\mathbf{Z}^{(*)})=0$ and
	\begin{equation}\label{Th3}
	\begin{aligned}
	\lim_{t\rightarrow \infty	}\| \nabla f(\mathbf{J}^{(t)},\mathbf{Z}^{(t)},\tilde{\mathbf{X}}^{(m)(t)})\|&=0.
	\end{aligned}
	\end{equation}
\end{theorem}

\begin{proof}
	First we prove that the sequence
	$\{(\mathbf{J}^{(t)},\mathbf{Z}^{(t)},\tilde{\mathbf{X}}^{(m)(t)})\}$ converges. Because  $(\mathbf{J}^{(t)},\mathbf{Z}^{(t)},\tilde{\mathbf{X}}^{(m)(t)})$ is bounded, there exists at least one accumulation point $(\mathbf{J}^{(*)},\mathbf{Z}^{(*)},\tilde{\mathbf{X}}^{(m)(*)})$ such that a sub-sequence $\{\mathbf{J}^{(t_j)},\mathbf{Z}^{(t_j)},\tilde{\mathbf{X}}^{(m)(t_j)}\}$ converges to $(\mathbf{J}^{(*)},\mathbf{Z}^{(*)},\tilde{\mathbf{X}}^{(m)(*)})$, i.e., $$\mathbf{J}^{(t_j)} \rightarrow \mathbf{J}^{(*)}, \mathbf{Z}^{(t_j)} \rightarrow \mathbf{Z}^{(*)}, \tilde{\mathbf{X}}^{(m)(t_j)}\rightarrow \tilde{\mathbf{X}}^{(m)(*)}.$$
	On the other hand, we have $$\mathbf{J}^{(t+1)}-\mathbf{J}^{(t)}\rightarrow 0,$$ $$\mathbf{Z}^{(t+1)}-\mathbf{Z}^{(t)}\rightarrow 0,$$ $$\tilde{\mathbf{X}}^{(m)(t+1)}-\tilde{\mathbf{X}}^{(m)(t)}\rightarrow 0$$ from Lemma \ref{lemma2}. Thus, the sequence
	$\{(\mathbf{J}^{(t)},\mathbf{Z}^{(t)},\tilde{\mathbf{X}}^{(m)(t)})\}$ converges and (Eq.\eqref{Th1}) holds. From Theorem 4.4 of \cite{cai2010singular}, the solution $\mathbf{Z}^{(t)}$ of sub-problem (Eq.\eqref{eq28}) converges to $\mathbf{Z}^{(*)}$ satisfying $diag(\mathbf{Z}^{(*)})=0$.
	
	Next we show that the partial derivative sequences with respect to $\mathbf{J}$, $\mathbf{Z}$ and $\tilde{\mathbf{X}}^{(m)}$ converge to zero. It can be deduced from the solution process of the sub-problems that
	\begin{equation}\label{eq26}
	\begin{aligned}
	\nabla f_\mathbf{J}(\mathbf{J}^{(t+1)},\mathbf{Z}^{(t)},\tilde{\mathbf{X}}^{(m)(t)})&=0,\\
	\nabla f_{\tilde{\mathbf{X}}^{(m)}}(\mathbf{J}^{(t+1)},\mathbf{Z}^{(t)},\tilde{\mathbf{X}}^{(m)(t+1)})&=0.\\
	\nabla f_\mathbf{Z}(\mathbf{J}^{(t+1)},\mathbf{Z}^{(t+1)},\tilde{\mathbf{X}}^{(m)(t+1)})&=0.
	\end{aligned}
	\end{equation}
	Since the sequence  $\{ \mathbf{J}^{(t)}; \mathbf{Z}^{(t)};\tilde{\mathbf{X}}^{(m)(t)} \}$ converges, the conclusion (Eq. \eqref{Th3}) is obtained immediately.
\end{proof}

\textbf{Time complexity}: We analyze the time complexity by counting the number of elementary operations performed in each sub-problem. (1)  For the sub-problem of $\mathbf{J}$, the time and space-consuming Kronecker product in (Eq.(\ref{4.1-2})) is computed via the product form in (Eq.(\ref{4.1-1})) during the iteration process in the conjugate gradient algorithm. Since we employ the preconditioned conjugate gradient method, the algorithm usually terminates quickly. Its computational complexity is $O(n^3+n^2 d_i)$, where $d_i$ represents the feature numbers of the $i$th view.  (2) The solution of $\mathbf{X}_i$ has a closed form whose computational complexity is $O(n_m^3+d_inn_m)$, where $n_m$ represents the number of unobservable entries in the $i$th view. (3) Because the minimization problem in (Eq.(\ref{Uz})) is solved analytically, the Uzawa’s algorithm for computing $\mathbf{Z}$ in the sub-problem (Eq.(\ref{eq28})) only costs few iterations in practice. The computational complexity is $O(n^3)$ for obtaining $\mathbf{Z}.$ (4) The computational complexity of spectral clustering in the last step is $O(n^3).$  To sum up, the total time complexity of the AWSR method is $O(n^3)$. \textbf{Space complexity}: The matrices involved are the $i$th view data $\mathbf{X}_i\in \mathbb{R}^{d_i\times n},$ the self-representation matrix $\mathbf{Z}\in \mathbb{R}^{n\times n},$ the weighting matrix $\mathbf{W}^{(i)}\in \mathbb{R}^{d_i\times n}$ and the separated variable $\mathbf{J}\in \mathbb{R}^{n\times n}$ for $i=1,\ldots,v.$ Thus the space complexity is $O(n^2).$

\section{Experiment}
In this section, we demonstrate the performance of our AWSR method compared to performance of other widely-used methods on five datasets.


\subsection{Datasets}
The proposed AWSR method as well as other approaches is implemented on five benchmark datasets, named \text{ORL}\footnotemark[1] \footnotetext[1]{http://cam-orl.co.uk/facedatabase.html}, \text{Still}\footnotemark[2] \footnotetext[2]{https://www.di.ens.fr/willow/research/stillactions/}, \text{BBCSport}\footnotemark[3] \footnotetext[3]{http://mlg.ucd.ie/datasets/bbc.html}, \text{Olympics}\footnotemark[4] \footnotetext[4]{http://mlg.ucd.ie/aggregation/}, \text{Leaves}\footnotemark[5] \footnotetext[5]{https://github.com/cswanghao/gbs/blob/}. For the process of creating missing view, we refer to the processing method in literature \cite{liu2021self}. First, the $n_m=[n*r]$ data entries of the first view are deleted, where $n$ is the total number of entries and $r$ is the missing ratio. Then the $n_m$ data entries of the second view are deleted and those existing in the first view are removed to ensure that all data is observed at least in one view. Finally, $n_m$ data entries arbitrarily chosen from the other views are deleted. In the numerical experiments, we consider the missing ratio as 0.1, 0.2, 0.3, 0.4, 0.5. The details of each dataset are as follows.

\text{ORL:} The dataset  consists of 400 images about 40 different people at different times, lighting, facial expressions (eyes open, eyes closed, smiling and not smiling) and facial details (with/without glasses). Refer to \cite{samaria1994parameterisation} and  \cite{zhang2015low}, we use part of the images to generate the datasets \text{ORL1} and \text{ORL2}. \text{ORL1} is a dataset of dual views with dimensions 1024 and 288, respectively. Data in \text{ORL2} has three features, whose dimensions are 4096, 3304 and 6750, respectively.

\text{Still:} The dataset \cite{delaitre2010recognizing} consists of 467 images taken of six different movements, including grabbing, running, walking, throwing, squatting and kicking. In the experiment, \text{Still} is a three-view dataset with dimensions 200, 200 and 200 respectively.

\text{BBCSport:} The dataset \cite{greene2006practical} consists of 737 sports articles published on the BBC Sport website by journalists in five different areas including athletics, cricket, football, rugby and tennis. In the experiment, we only selecte a subset of \text{BBCSport} that contains news reports from 116 journalists. \text{BBCSport} is a four-view dataset with dimensions 1991, 2063, 2113 and 2158.

\text{Olympics:} The dataset \cite{greene2013producing} contains images of 464 athletes or organizers at 28 different sports in the London 2012 Summer Olympic Games. In the experiment, \text{Olympics} is a double\cl{-}view dataset with dimensions 4942 and 3097 respectively.

\text{Leaves:} This dataset \cite{beghin2010shape} contains 1600 leaf images, covering the leaf images of 100 plant species. \text{Leaves} describes leaves from shape, fine-scale edges and texture histogram features. 
In the experiment, \text{Leaves} is a three-view dataset with dimensions 64, 64 and 64.

Statistical information about the data is presented in Table \ref{tab2}. The  ``Number of Clusters" and ``Dimension of Each View'' are acquired from the information of the dataset. For example, \text{ORL} dataset is the data of 400 images composed of different images of 40 people, and its number of clusters is 40.
\begin{table}[htb]
		\caption{ Statistics of the datasets.}\label{tab2}	
	\begin{tabular}{@{}llllll@{}}
		\toprule
		\multirow{2}{*}{Dataset} & Number  of & \multicolumn{4}{c}{Dimension of Each View} \\
		  \cmidrule{3-6}
  & Clusters & 1& 2&3 &4\\
 \midrule
		\text{ORL1}                        & 40     & 1024  & 288   & -    & -    \\
		\text{Still}                      & 6      & 200   & 200   & 200  & -    \\
		\text{BBCSport}                   & 5      & 1991  & 2063  & 2113 & 2158 \\
		\text{Olympics}                   & 28     & 4942  & 3097  & -    & -    \\
		\text{Leaves}                    & 100    & 64    & 64    & 64   & -    \\
		\text{ORL2}              &40      & 4096  & 3304  & 6750 & -   \\
		 \bottomrule
	\end{tabular}
\end{table}

\subsection{Methods for comparison}
In this subsection, we provide a brief introduction to eight other algorithms and two benchmarks, which we  use to
with our proposed algorithm.

(1) LSRs (single-best baseline) \cite{liu2012robust} first filled the missing data with random values, then it applied LSR algorithm to the learning of each view.

(2) In LSRc (concatenated baseline) \cite{liu2012robust}, the missing data was filled with arbitrary values, and the data of each view was
concatenated. It run the LSR algorithm on the concatenated data.

(3) IMG \cite{zhao2016incomplete} factorized data of each view to learn the latent subspaces independently. Meanwhile, it utilized the graph Laplacian term to regularize the latent subspaces of each view.

(4) DAIMC \cite{hu2019doubly} was a weighted non-negative matrix factorization based method. It employed the existing instance alignment information to learn the consistent latent feature matrix of all views. Besides, in order to reduce the influence of missing entries,
the $\ell_{2,1}$-norm was applied to build the consistency basis matrix of all views.

(5) AGL \cite{wen2018incomplete} first performed a low rank representation of the graph for each view with spectral constraints. Then it built a consensus representation for all views using the co-regularization term.

(6) AWGF \cite{zhang2020adaptive} factorized a weighted non-negative matrix to obtain the feature matrix of each view and 
then constructed the graph of each view. The feature extraction and the graph were fused into a large framework according to a certain weight.

(7) PLR \cite{lian2021partial} obtained the consensus feature matrix of all views based on non-negative matrix factorization. 
Then it imposed $\ell_{2,1}$-norm constraints on the basis matrix of each view and finally imposed regularization constraints on the local graph to obtain the shared feature matrix.

(8) IMSR \cite{liu2021self} performed feature extraction. 
It employed missing data imputation and low rank regularized consensus coefficient matrix via self-representation learning to obtain the complete consensus coefficient matrix.

(9) CAMVSC \cite{tang2022consistent} employed an auto-weighted technique  and constructed a self-representation subspace model   to cluster mult-views on complete dataset.

(10) SIMC\_ADC \cite{he2023scalable} proposed an anchor point learning method to obtain consensus coefficient matrix for  large scale multi-view clustering problem.


\subsection{Evaluation metrics}
Accuracy (ACC), normalized mutual information (NMI) and Purity are used to measure the performance of clustering. Generally speaking, higher values represent better clustering performance. Next, we introduce the definitions of the metrics. The ACC is calculated by
\begin{equation*}
\begin{aligned}
ACC = \frac{TP+TN}{TP+TN+FP+FN} ,
\end{aligned}
\end{equation*}
where true positive (TP) indicates that similar documents are put in the same cluster, true negative (TN) means that different documents are placed in various clusters, false negative (FN) indicates that similar documents are placed in different clusters and false positive (FP) indicates that various documents are put in the same cluster. The mutual information (MI) is
\begin{equation*}
\begin{aligned}
MI (\mathcal{C},\varTheta) = \sum\limits_{c_i\in \mathcal{C}} \sum\limits_{\omega_k\in \mathcal{\varTheta}} P(c_i,\omega_k)log \frac{P(c_i,\omega_k)}{P(c_i)P(\omega_k)} ,
\end{aligned}
\end{equation*}
where $P(\cdot)$ is a probability function. If $\mathcal{C}$ and $\varTheta$ are independent of each other, mutual information (MI) is equivalent to zero. In the clustering experiments, we generally use the normalized mutual information (NMI) to measure
\begin{equation*}
\begin{aligned}
NMI (\mathcal{C},\varTheta) = \frac{MI(\mathcal{C},\varTheta)}{max(H(\mathcal{C}),H(\varTheta))} ,
\end{aligned}
\end{equation*}
where $H(\cdot)$ is the entropy function and 
defined as $H(\varTheta)=-\sum\limits_{k}P(\omega_k)logP(\omega_k)$.
The metric Purity is given by
\begin{equation*}
\begin{aligned}
Purity (\mathcal{C},\varTheta) = \frac{1}{n}\sum\limits_{i=1}^{r} \max\limits_{j=1}^k (\omega_k \cap c_i),
\end{aligned}
\end{equation*}
where $n$ is the total number of sample, $\varTheta=\{\omega_1,\omega_2,\cdot \cdot \cdot,\omega_K\}$ is the clustering class and the real class is defined as $\mathcal{C} = \{c_1,c_2,\cdot \cdot \cdot,c_J\}$.

\subsection{Numerical experiments}

In this section, we explain how to choose the parameters. Then we demonstrate the numerical results of the proposed AWSR method comparing with other prevailing methods. At last, the sequence of  objective function values is plotted to verify the convergence property of the computational algorithm.

\subsubsection {Parameter selection}
For the elements $\alpha_{mj}^i$ in the weighting matrix, we assign the constant $1$ to it from experience. The main idea we determine the parameters in the objective function is that we compare the performance of the AWSR method with different potential candidates and choose the best ones as the value of parameters.

First, we list the potential values of the parameters. We roughly approximate the range  which may contain the desirable value of the  parameters. All potential values of the parameters are expressed by $2^q.$  In Table \ref{para_data}, we demonstrate the values of the exponential  $q$ for each parameter in different datasets. To be specific, the values assigned for the parameter $\gamma$ in the \text{ORL1} dataset are $2^{3.1}, 2^{3.3},2^{3.5}.$  That is because the expression $[3.1:0.2:3.5]$ indicates that the exponential $q$ takes the values from 3.1 to 3.5 with a common difference 0.2 between any successive values. Other nominated parameters are chosen based on Table \ref{para_data} in the same way.

\begin{table}[htpb]
	\setlength\tabcolsep{0.05pt}
	\caption{The values of the exponential $q$ for different parameters ($2^q$).\label{para_data}}
	\begin{tabular}{@{}l|l|l|l@{}}
		\toprule
        Dataset & $ q $ for $\lambda=2^q $        &  $ q$ for $ \gamma=2^q $      &   $ q $ for $\alpha=2^q $\\
\hline
		\text{ORL1}	& $[-2.9:0.1:-1.7]$ & $[3.1:0.2:3.5]$ & $[7.3:0.1:8]$\\
	\hline	
		\text{Still} & $[-1.9:0.05:-1.7]$ & $[3.1:0.2:3.5]$ & $[7.2:0.1:7.7]$\\
\hline	
		\text{BBCSport}& $[-10:2:10]$ & $[2.5:1:4.5]$ & $[7.3:0.1:8]$\\
	\hline	
		\text{Olympics}& $[-1.8:0.1:-1.7]$ & $[3.1:0.2:3.5]$ & $[7.2:0.1:8]$\\
	\hline	
		\text{Leaves}& $[-2.3:0.1:-1.6]$ & $[3.1:0.2:3.5]$ & $[7.4:0.1:8]$\\
	\hline
       \text{ORL2} & [-3:1:0] & [2:1:6] & [6.7:0.05:6.85]\\
\midrule
	\end{tabular}
\end{table}	

\begin{figure}[htpb]
	\centering
	\includegraphics[width=0.5\textwidth]{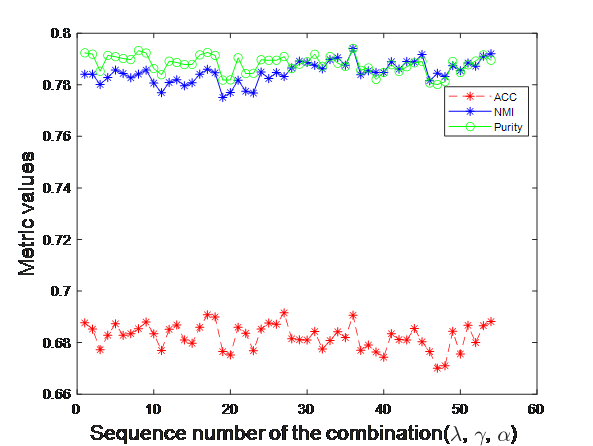}
	\caption{Comparison of numerical results of different parameter combinations on Olympics.}\label{Olymp_para}
\end{figure}

Next, we take the Olympics dataset as an example to show the process of choosing parameters. As presented in the `Olympics' row in Table \ref{para_data}, the parameter $\lambda$ has 2 candidates, $2^{-1.8}$ and $2^{-1.7},$ while $\gamma$ is allowed to take the values of $2^{3.1}, 2^{3.3},2^{3.5}.$ The parameter $\alpha$ has 9 choices, from $2^{7.2}$ to $2^8.$ Thus the combination of $(\lambda,\gamma,\alpha)$ has 54 cases. We sort these 54 combinations in lexical order, which means every number from 1 to 54  corresponds to a combination of $(\lambda,\gamma,\alpha)$. For each combination of $(\lambda,\gamma,\alpha),$ we run the AWSR method when the missing ratio is 0.1, 0.2 , 0.3 and 0.4 respectively. In Figure \ref{Olymp_para}, the average values of  ACC, NMI, Purity are demonstrated. The number in the horizontal axis is the order of the corresponding combination   $(\lambda,\gamma,\alpha)$.  It can be seen that, the model performs the best when  the parameters come to the $36$th combination of  $(\lambda,\gamma,\alpha),$ which is $\lambda=2^{(-1.7)}$, $\gamma=2^{(3.1)}$ and $\alpha=2^{(8.0)}$.

\subsubsection {Numerical results}

\textbf{Complete datasets}: First, we compare the performance of our proposed algorithm with CAMVSC algorithm on complete datasets \text{ORL2,} \text{BBCSport} and \text{Leaves}. The experimental results are shown in Table \ref{tab4}. Since the CAMVSC algorithm is designed especially for complete multi-view clustering, it outperforms the AWSR method on these complete dataset.
\begin{table}[htpb]
	\caption{ Performance comparison between the algorithm proposed and CAMVSC algorithm on complete dataset \text{ORL2}, \text{BBCSport} and \text{Leaves}.} \label{tab4}
	\begin{tabular}{@{}lllll@{}}
		\toprule
		\multicolumn{1}{l}{Dataset}    & algorithm  & \multicolumn{1}{c}{ACC}  & NMI  & Purity  \\
		\midrule
		\multirow{1}{*}{\text{ORL2}}    & AWSR      & \multicolumn{1}{c}{80.90} & 91.16 & 84.30 \\
		& CAMVSC    & 82.45                     & 91.73 & 85.33\\ \midrule
		\multirow{1}{*}{\text{BBCSport}}    & AWSR      & 80.17                     & 75.89 & 90.52\\
		& CAMVSC    & 92.24                      & 82.19  & 92.24\\ \midrule
		\multirow{1}{*}{\text{Leaves}} & AWSR      & 77.80                     & 89.47 & 80.30  \\
		& CAMVSC    &91.38                     & 81.05 & 91.38 \\
		\midrule
	\end{tabular}
\end{table}

\begin{figure*}[htbp]	
	\begin{minipage}{0.19\linewidth}
		\centerline{\text{Leaves}}
		\vspace{3pt}
		\centerline{\includegraphics[width=\textwidth]{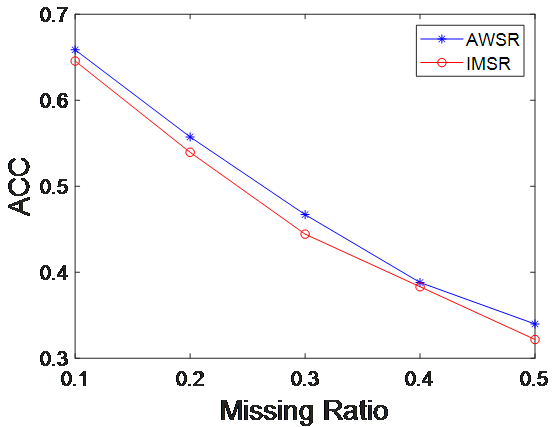}}
		\vspace{3pt}
		\centerline{\includegraphics[width=\textwidth]{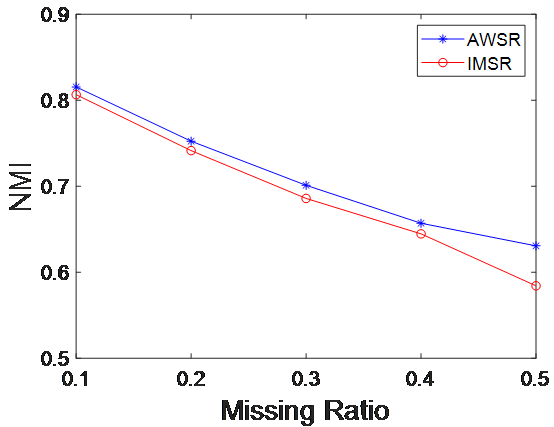}}
		\vspace{3pt}
		\centerline{\includegraphics[width=\textwidth]{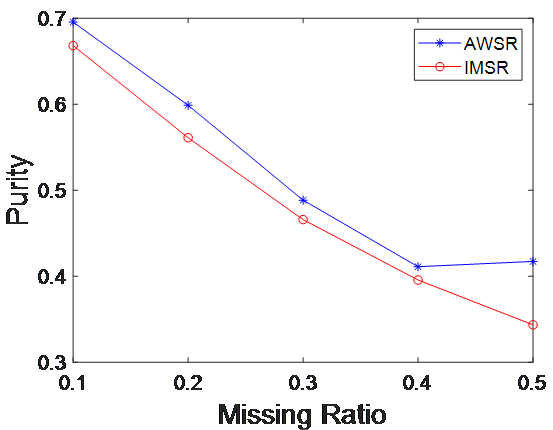}}
		\vspace{3pt}	
	\end{minipage}	
	\begin{minipage}{0.19\linewidth}
		\centerline{\text{BBCSport}}
		\vspace{2pt}
		\centerline{\includegraphics[width=\textwidth]{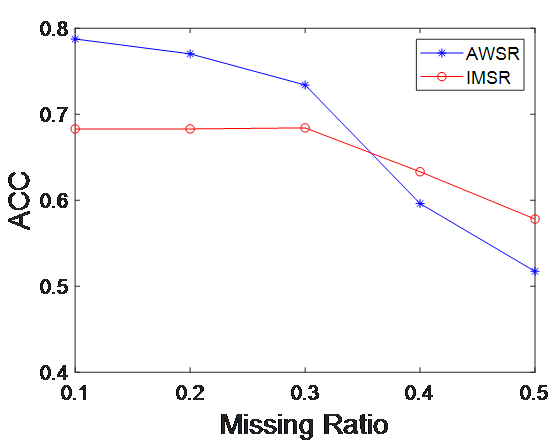}}
		\vspace{3pt}
		\centerline{\includegraphics[width=\textwidth]{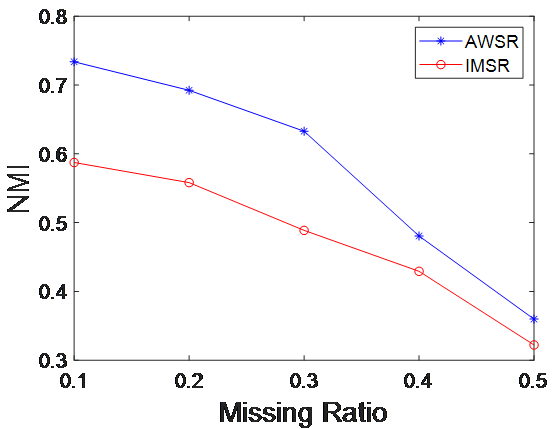}}
		\vspace{3pt}
		\centerline{\includegraphics[width=\textwidth]{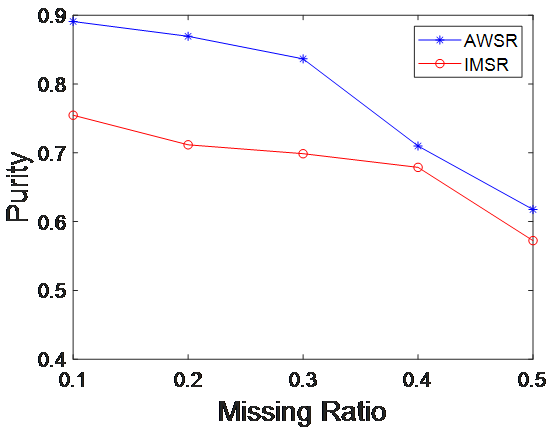}}
		\vspace{3pt}			
	\end{minipage}
	\begin{minipage}{0.19\linewidth}
		\centerline{\text{Olympics}}
		\vspace{3pt}
		\centerline{\includegraphics[width=\textwidth]{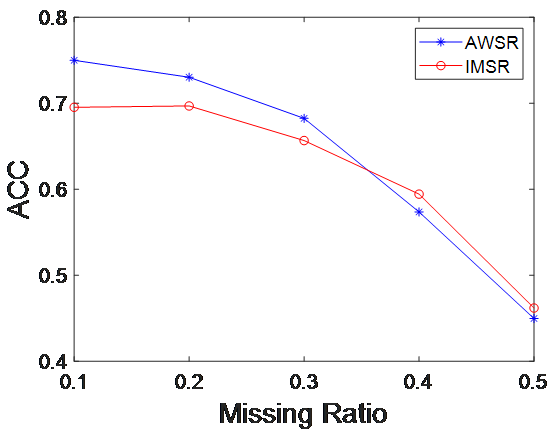}}
		\vspace{3pt}
		\centerline{\includegraphics[width=\textwidth]{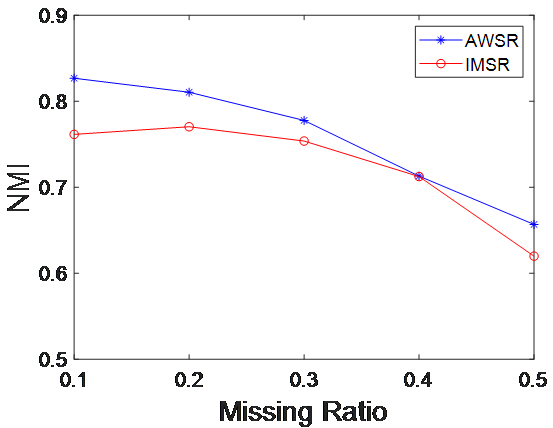}}
		\vspace{3pt}
		\centerline{\includegraphics[width=\textwidth]{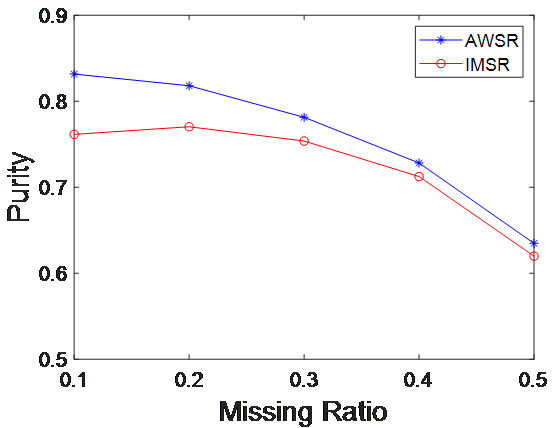}}
		\vspace{3pt}		
	\end{minipage}
    \begin{minipage}{0.19\linewidth}
        \centerline{\text{ORL1}}
		\vspace{3pt}
		\centerline{\includegraphics[width=\textwidth]{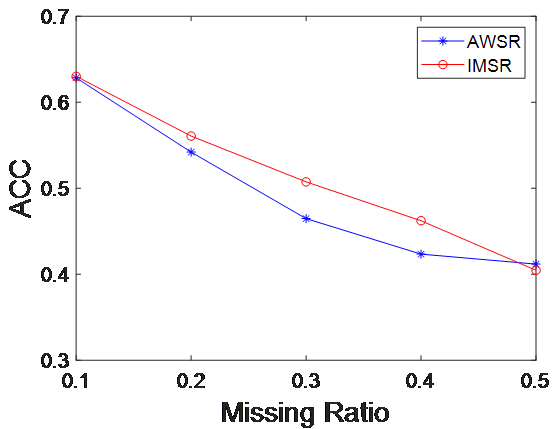}}
		\vspace{3pt}
		\centerline{\includegraphics[width=\textwidth]{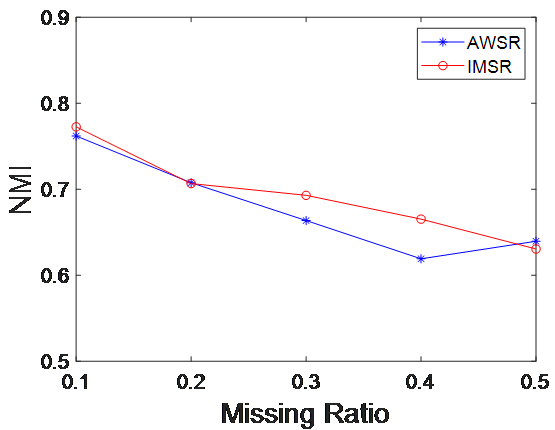}}
		\vspace{3pt}
		\centerline{\includegraphics[width=\textwidth]{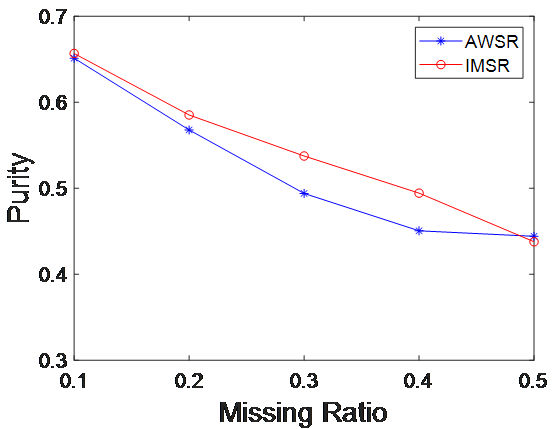}}
		\vspace{3pt}
    \end{minipage}
 \begin{minipage}{0.19\linewidth}
        \centerline{\text{Still}}
		\vspace{3pt}
		\centerline{\includegraphics[width=\textwidth]{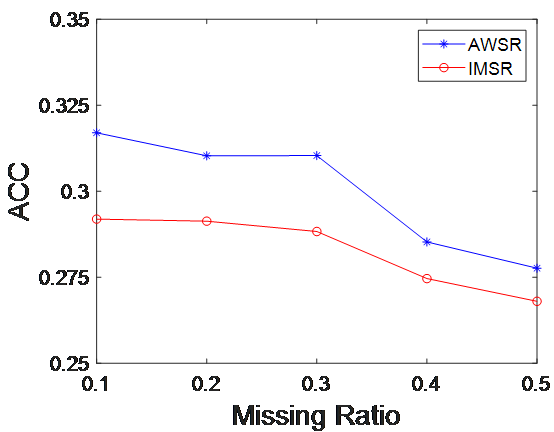}}
		\vspace{3pt}
		\centerline{\includegraphics[width=\textwidth]{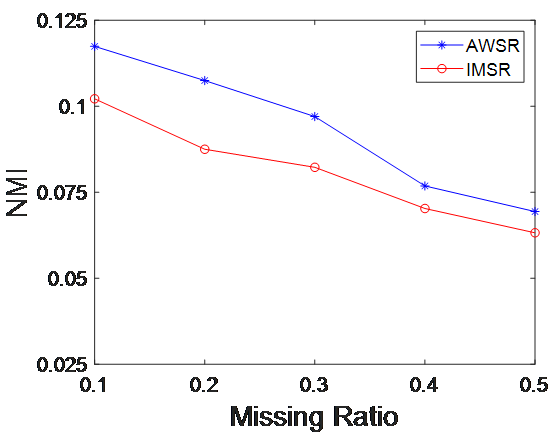}}
		\vspace{3pt}
		\centerline{\includegraphics[width=\textwidth]{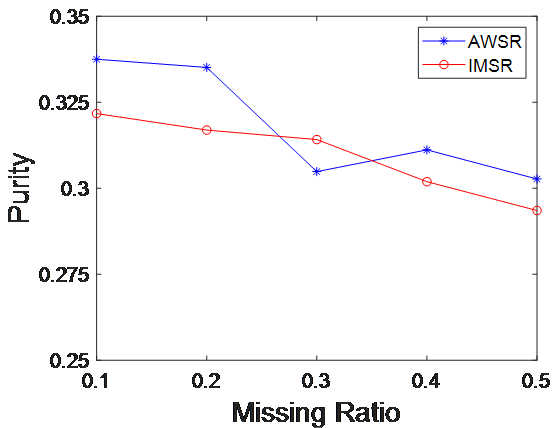}}
		\vspace{3pt}
    \end{minipage}
	\caption{Performance comparison of AWSR and IMSR methods on datasets \text{Leaves}, \text{BBCSport}, \text{Olympics}, \text{ORL1} and \text{Still}.}\label{fig3}
\end{figure*}

\noindent \textbf{ Incomplete datasets}: In \cite{liu2021self}, it is shown that the IMSR algorithm is superior to the LSRs, LSRc, IMG, DAIMC, AGL, AWGF and PLR for incomplete dataset. Therefore first we only compare the AWSR algorithm  with the IMSR algorithm  and  the state of the art SIMC\_ADC algorithm when the missing ratio changes from 0.1 to 0.5. The results are shown in Figure \ref{fig3} and Table \ref{SIAW} respectively. In Figure \ref{fig3}, the values of ACC, NMI and Purity computed by AWSR and IMSR are plotted. The blue lines describe the evaluation indices given by AWSR, while the red ones show the evaluation indices generated by IMSR. In Table \ref{SIAW}, the superior values are in bold type. For most cases,  AWSR performs better than IMSR and SIMC\_ADC.

\begin{figure*}[h]	
	\begin{minipage}{0.33\linewidth}
\centerline{\text{Leaves}}
		\vspace{3pt}
		\centerline{\includegraphics[width=\textwidth]{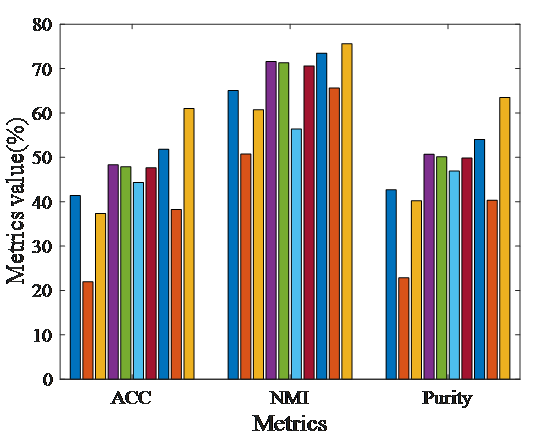}}
		\centerline{\text{Still}}
		\vspace{3pt}
		\centerline{\includegraphics[width=\textwidth]{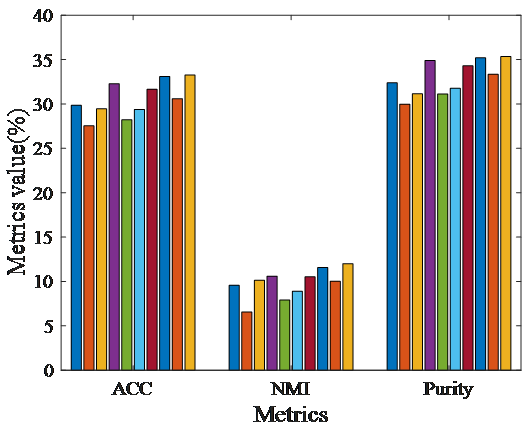}}
	\end{minipage}	
	\begin{minipage}{0.33\linewidth}
\centerline{\text{Olympics}}
		\vspace{3pt}
		\centerline{\includegraphics[width=\textwidth]{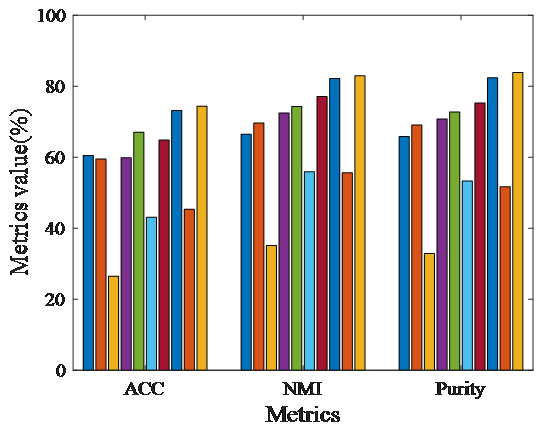}}
		\centerline{\text{ORL1}}
		\vspace{2pt}
		\centerline{\includegraphics[width=\textwidth]{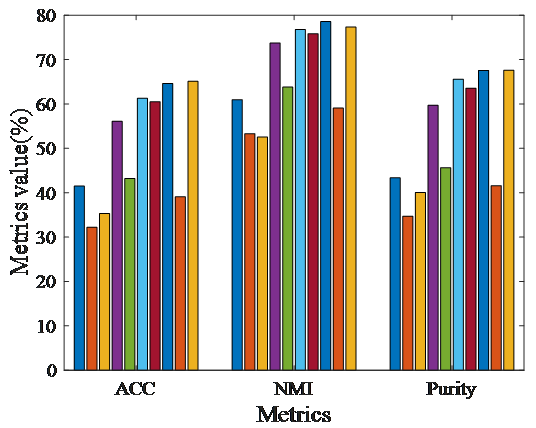}}		
	\end{minipage}
\begin{minipage}{0.33\linewidth}
		\centerline{\text{BBCSport}}
		\vspace{3pt}
		\centerline{\includegraphics[width=\textwidth]{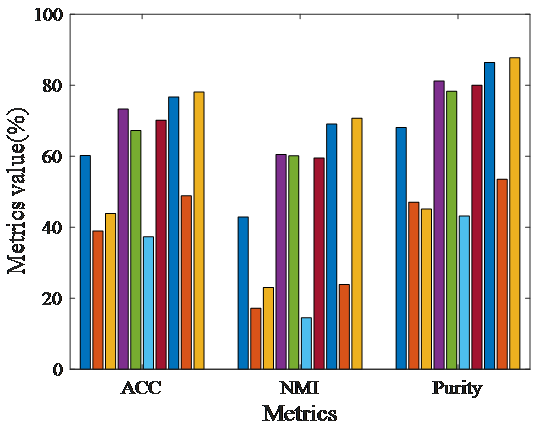}}
		\vspace{3pt}
		\centerline{\includegraphics[scale=0.62]{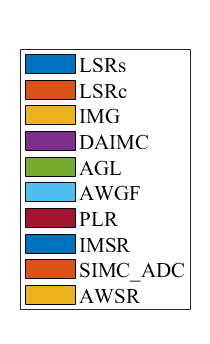}}
	\end{minipage}
	\caption{Visualization histogram of the results in Table \ref{tab3}.\label{metrics}}
\end{figure*}

	\begin{table*}[htbp]
\caption{Comparison between AWSR method and the SIMC\_ADC method on incomplete datasets with various missing ratios. The best results are shown in bold. \label{SIAW}}
		\setlength\tabcolsep{10pt}
		\begin{tabular}{llcccccc}
			\toprule
			 Metrics & Dataset &	Algorithm  & 0.1 & 0.2 & 0.3 & 0.4 & 0.5 \\
	\midrule
	\multirow{10}{*}{ \textbf{ACC}} &\multirow{2}{*}{ORL1} & SIMC\_ADC & 51.73 & 45.85& 37.80& 29.23&  30.80   \\
	\ & \  &	AWSR & \textbf{65.38} &\textbf{55.15} & \textbf{47.90} & \textbf{43.73} & \textbf{41.98}\\
\cmidrule{2-8}
	\ &  \multirow{2}{*}{ Olympics} &	SIMC\_ADC &56.51  &51.44  &45.60  &39.01  &34.24\\
	\ &	\ &	AWSR &\textbf{76.57}  & \textbf{73.69} &\textbf{68.32}  & \textbf{57.74} &\textbf{44.29} \\
		\cmidrule{2-8}
\ &  \multirow{2}{*}{ Leaves} &	SIMC\_ADC &52.62  &42.84  &36.36  &32.03  &27.36\\
	\ &	\ &	AWSR &\textbf{66.25}  &\textbf{56.48} &\textbf{46.31}  & \textbf{39.71} &\textbf{33.99} \\
\cmidrule{2-8}
\ &  \multirow{2}{*}{ BBCSport} &	SIMC\_ADC &52.33  &49.31  &48.97  &48.53  &45.07\\
	\ &	\ &	AWSR &\textbf{79.31}  &\textbf{76.55} &\textbf{71.53} &\textbf{61.46}  & \textbf{53.97} \\
\cmidrule{2-8}
\ &  \multirow{2}{*}{ Still} &	SIMC\_ADC &32.63  &\textbf{31.32}  &29.13  &\textbf{29.79} &\textbf{29.14}\\
	\ &	\ &	AWSR &\textbf{33.02}  &31.18 &\textbf{29.87}  &28.56 &27.20 \\
		\cmidrule{1-8}
	\multirow{10}{*}{ \textbf{NMI}} &\multirow{2}{*}{ORL1} & SIMC\_ADC & 69.64 & 65.20& 57.63& 50.17&  52.90   \\
	\ & \  &	AWSR & \textbf{78.01} &\textbf{72.43} & \textbf{67.13} & \textbf{62.94} & \textbf{64.22}\\
\cmidrule{2-8}
	\ &  \multirow{2}{*}{ Olympics} &	SIMC\_ADC &65.84  &61.25  &55.27  &49.80  &45.86\\
	\ &	\ &	AWSR &\textbf{83.35}  & \textbf{81.36} &\textbf{78.19}  & \textbf{71.35} &\textbf{65.33} \\
		\cmidrule{2-8}
\ &  \multirow{2}{*}{ Leaves} &	SIMC\_ADC &74.71  &69.21  &64.70  &61.45  &58.07\\
	\ &	\ &	AWSR &\textbf{81.48}  &\textbf{75.47} &\textbf{69.72}  & \textbf{66.00} &\textbf{60.00} \\
\cmidrule{2-8}
	\ &  \multirow{2}{*}{ BBCSport} &	SIMC\_ADC &28.65  &23.47  &23.71  &22.88  &20.61\\
	\ &	\ &	AWSR &\textbf{73.60}  &\textbf{68.46} &\textbf{63.41}  & \textbf{50.11} &\textbf{39.41} \\
\cmidrule{2-8}
\ &  \multirow{2}{*}{ Still} &	SIMC\_ADC &11.93  &\textbf{10.54}  &\textbf{9.74}  &\textbf{9.21} &\textbf{8.70}\\
	\ &	\ &	AWSR &\textbf{12.08}  &10.25 &8.97  & \textbf{}7.80 &6.67 \\	
\cmidrule{1-8}
	\multirow{10}{*}{ \textbf{Purity}} &\multirow{2}{*}{ORL1} & SIMC\_ADC & 54.39 & 48.78& 40.48& 31.38&  32.95   \\
	\ & \  &	AWSR & \textbf{68.25} &\textbf{57.70} & \textbf{50.73} & \textbf{46.48} & \textbf{45.18}\\
\cmidrule{2-8}
	\ &  \multirow{2}{*}{ Olympics} &	SIMC\_ADC &64.38  &58.92  &52.97  &39.09  &43.82\\
	\ &	\ &	AWSR &\textbf{83.95}  & \textbf{82.18} &\textbf{78.71}  & \textbf{72.72} &\textbf{62.85} \\
		\cmidrule{2-8}
\ &  \multirow{2}{*}{ Leaves} &	SIMC\_ADC &54.73  &45.10  &38.68  &39.94  &29.16\\
	\ &	\ &	AWSR &\textbf{68.27}  &\textbf{58.63} &\textbf{48.20}  & \textbf{41.96} &\textbf{39.79} \\
\cmidrule{2-8}
\ &  \multirow{2}{*}{ BBCSport} &	SIMC\_ADC & 58.71 & 55.35 &52.67  &52.50  &48.19\\
	\ &	\ &	AWSR &\textbf{89.40}  &\textbf{86.38} &\textbf{84.40}  & \textbf{72.07} &\textbf{64.40} \\
\cmidrule{2-8}
\ &  \multirow{2}{*}{ Still} &	SIMC\_ADC &\textbf{35.56}  &\textbf{34.28} &\textbf{33.00}  & \textbf{32.35} &\textbf{31.50} \\
	\ &	\ &	AWSR 		 & 35.25 & 33.43 &32.29  &30.90 & 29.27\\
\bottomrule
\end{tabular}

\end{table*}
Moreover, we also compare the AWSR method with the above-mentioned  approaches. The average values of the clustering results provided with missing ratio being 0.1, 0.2, 0.3, 0.4 and 0.5 are demonstrated in Table \ref{tab3}.  The best metric value of each problem is shown in bold. We can see that for each dataset, the overall performance of the proposed algorithm AWSR on ACC, NMI and Purity is significantly better than that of other algorithms except for the NMI metric on \text{ORL1} dataset, of which the proposed algorithm is the second best.
\begin{table*}[htbp]
\setlength\tabcolsep{4pt}
    \caption{ Comparison of the performance after collation. The average values of the clustering results provided with missing ratio being 0.1, 0.2, 0.3, 0.4 and 0.5 are given.  The best results are shown in bold.}\label{tab3}
	\begin{tabular}{@{}llllllllllcl@{}}
		\toprule
  	\multicolumn{1}{c}{}    & Dataset  & \multicolumn{1}{c}{LSRs}  & LSRc  & IMG   & DAIMC & AGL   & AWGF  & PLR   & IMSR &SIMC\_ADC & AWSR    \\
\midrule
		\multirow{5}{*}{ACC}    & \text{ORL1}      & \multicolumn{1}{c}{41.50} & 32.20 & 35.30 & 56.10 & 43.20 & 61.27 & 60.50 & 64.63 &39.08 &\textbf{65.10}  \\
		& \text{Still}    & 29.85                     & 27.54 & 29.46 & 32.29 & 28.22 & 29.38 & 31.65 & 33.09 & 30.60& \textbf{33.26}   \\
		& \text{BBCSport} & 60.17                     & 38.97 & 43.90 & 73.28 & 67.24 & 37.29 & 70.17 & 76.72 &48.84& \textbf{78.10}  \\
		& \text{Olympics} & 60.52                     & 59.53 & 26.47 & 59.87 & 67.03 & 43.09 & 64.87 & 73.17 &45.36& \textbf{74.37}  \\
		& \text{Leaves}   & 41.39                     & 21.96 & 37.36 & 48.33 & 47.88 & 44.34 & 47.63 & 51.81 & 38.24& \textbf{61.02}  \\ \midrule
		\multirow{5}{*}{NMI}    & \text{ORL1}      & 60.98                     & 53.31 & 52.55 & 73.73 & 63.80 & 76.79 & 75.83 & \textbf{78.57} & 59.08 & 77.35   \\
		&\text{Still}    & 9.56                      & 6.56  & 10.13 & 10.60 & 7.91  & 8.90  & 10.52 & 11.56 &10.02& \textbf{12.01}   \\
		& \text{BBCSport} & 42.88                     & 17.16 & 23.04 & 60.49 & 60.10 & 14.51 & 59.51 & 69.06 &23.86& \textbf{70.7}   \\
		& \text{Olympics} & 66.45                     & 69.67 & 35.12 & 72.44 & 74.28 & 55.92 & 77.09 & 82.16 &55.60& \textbf{82.94}  \\
		& \text{Leaves}   & 65.06                     & 50.75 & 60.75 & 71.59 & 71.32 & 56.40 & 70.61 & 73.49 &65.63& \textbf{75.63}  \\ \midrule
		\multirow{5}{*}{Purity} & \text{ORL1}      & 43.35                     & 34.70 & 40.06 & 59.70 & 45.65 & 65.58 & 63.55 & 67.57 & 41.59 &\textbf{67.61}  \\
		& \text{Still}    & 32.38                     & 29.98 & 31.15 & 34.90 & 31.13 & 31.76 & 34.30 & 35.20 &33.34& \textbf{35.35}   \\
		&\text{BBCSport} & 68.10                     & 47.07 & 45.10 & 81.21 & 78.28 & 43.19 & 80.00 & 86.38 &53.48& \textbf{87.75}  \\
		& \text{Olympics} & 65.82                     & 69.09 & 32.88 & 70.78 & 72.76 & 53.32 & 75.30 & 82.37 &51.68& \textbf{83.86}  \\
		&\text{Leaves}   & 42.69                     & 22.85 & 40.23 & 50.68 & 50.13 & 46.91 & 49.85 & 54.03 &40.32& \textbf{63.51}  \\
\midrule

\end{tabular}
\end{table*}
We also visualize the results of Table \ref{tab3} with different metrics values based on different datasets  in Figure \ref{metrics}.

\vspace{2pt}
\noindent \textbf{Time comparison}: In terms of computation efficiency, the AWSR method is not as efficient as the IMSR and  SIMC\_ADC algorithm. Take Olympics, BBCSport and Still as examples.  In Table \ref{time} we show the  the  computation time of  AWSR, IMSR  and  SIMC\_ADC   when the missing rate is $0.1$.
 \begin{table}[htbp]
\setlength\tabcolsep{2pt}
		\caption{Comparison of  time (seconds).}\label{time}	
	\begin{tabular}{@{}llll@{}}
		\toprule
		\diagbox{Method}{Time}{Dataset}  &\text{Olympics}&\text{BBCSport}&\text{Still} \\
	\midrule	
         AWSR & 6.596  &1.656&5.159\\
		IMSR         &2.166    & 0.378    & 0.724\\
		SIMC\_ADC    & 2.174      &1.097 & 1.439    \\
		 \bottomrule
	\end{tabular}
\end{table}

{\subsubsection {Sequence of objective function values}

To verify the convergence property of the BCD algorithm, we  show the result of objective function values on \text{Still}, \text{Leaves}, \text{ORL1} and \text{Olympics} datasets in Figure \ref{fig5}. Since the missing elements at the initial point are zero, it can be seen that the first objective function value in each figure is less than the objective function value at the second point, while some missing elements are nonzero. The objective function value declines and  tends to converge eventually.
\begin{figure}[h]	
	\begin{minipage}{0.49\linewidth}
		\centerline{\text{Still}}
		\vspace{3pt}
		\centerline{\includegraphics[width=\textwidth]{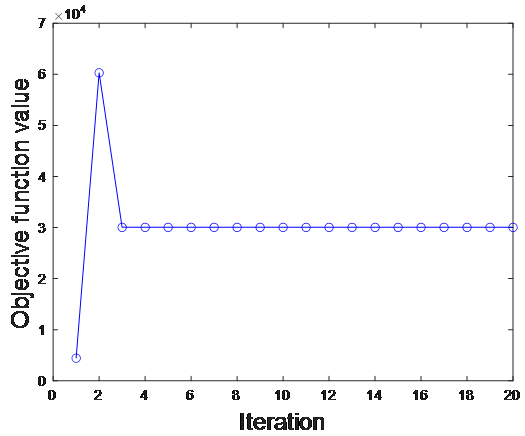}}
		\centerline{\text{Leaves}}
		\vspace{3pt}
		\centerline{\includegraphics[width=\textwidth]{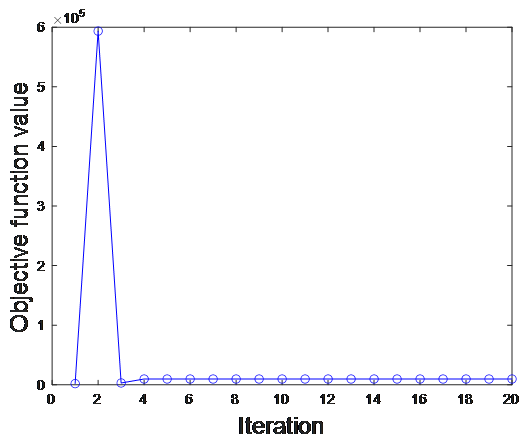}}
	\end{minipage}	
	\begin{minipage}{0.49\linewidth}
		\centerline{\text{ORL1}}
		\vspace{2pt}
		\centerline{\includegraphics[width=\textwidth]{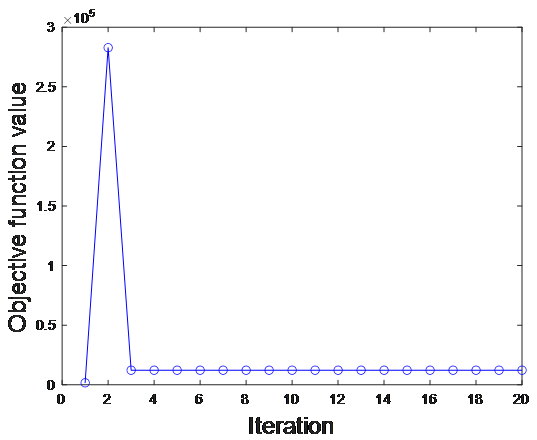}}
		\centerline{\text{Olympics}}
		\vspace{3pt}
		\centerline{\includegraphics[width=\textwidth]{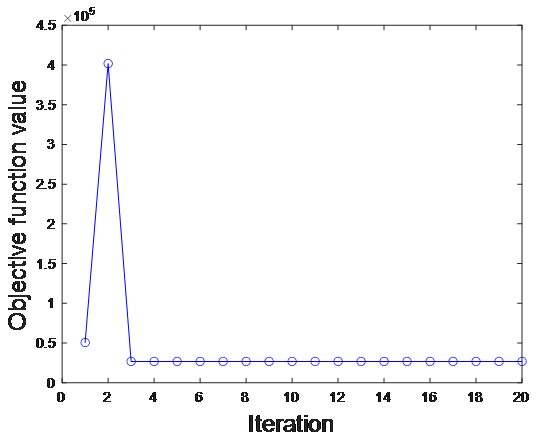}}
	\end{minipage}
	\caption{The convergence curves of the objective function value on datasets \text{Still}, \text{Leaves}, \text{ORL1} and \text{Olympics}.\label{fig5}}
\end{figure}

\section{Conclusion}
For problems of IMC, we introduce an adaptive weighted self-representation (AWSR) model. It adjusts the weighting matrix on the basis of the variations of views and the recovery process of the missing data. The proposed AWSR model is concise and can be solved by the traditional BCD algorithm. Convergence property of the iterative algorithm is proved. Numerical experiments on five real-world incomplete data demonstrate the effectiveness of AWSR and its superiority over other eight widely-used approaches. In addition, we also compare the performance with CAMVSC on three complete data, our method can be applied to complete and incomplete datasets, while the CAMVSC value can only be applied to complete ones.

Although the AWSR method is effective for IMC, there are still some issues that deserve further investigation.
Theoretically, the elements of weighting matrix are determined by the differences of instances, uncompleted entries and completed entries. However, the adjustment of the entries in the weighting matrix is realized empirically in the experiments. A more rational and flexible updating strategy should be provided. Besides, high-efficiency optimization algorithms are needed for large scale datasets.

\bmhead{Acknowledgments}
This work was supported by the National Natural Science Foundation of China [grant No.12326302, No.62073087].

\bmhead{Data availability}
The data used to support the findings of this study are available from the corresponding author upon request.


\section*{Declarations}
\begin{itemize}

\item \textbf{Conflict of interest} The authors declare that they have no known competing financial interests or personal relationships that could have appeared to influence the work reported in this paper.
\item \textbf{Author contributions}
 Conceptualization: [Lishan Feng], [Jingya Chang]; Software: [Lishan Feng]; Validation: [Lishan Feng], [Guoxu Zhou],  [Jingya Chang]; Formal analysis: [Lishan Feng]; Data curation: [Lishan Feng]; Writing - original draft: [Lishan Feng]; Preparation: [Lishan Feng]; Visualization: [Lishan Feng]; Project administration: [Lishan Feng]; Investigation: [Guoxu Zhou], [Jingya Chang]; Resources: [Guoxu Zhou],  [Jingya Chang]; Supervision:[Guoxu Zhou], [Jingya Chang]; Funding acquisition: [Guoxu Zhou], [Jingya Chang]; Methodology: [Jingya Chang]; Writing - review $\&$ editing: [Jingya Chang].
All authors have read and agreed to the published version of the manuscript.
\end{itemize}


\bibliographystyle{sn-basic}
\bibliography{ref}

\begin{thebibliography}{45}
\providecommand{\natexlab}[1]{#1}
\providecommand{\url}[1]{{#1}}
\providecommand{\urlprefix}{URL }
\providecommand{\doi}[1]{\url{https://doi.org/#1}}
\providecommand{\eprint}[2][]{\url{#2}}
 \bibcommenthead

\bibitem[{Beghin et~al(2010)Beghin, Cope, Remagnino, and
  Barman}]{beghin2010shape}
Beghin T, Cope JS, Remagnino P, et~al (2010) Shape and texture based plant leaf
  classification. In: International conference on advanced concepts for
  intelligent vision systems, Springer, pp 345--353

\bibitem[{Bramble et~al(1997)Bramble, Pasciak, and
  Vassilev}]{bramble1997analysis}
Bramble JH, Pasciak JE, Vassilev AT (1997) Analysis of the inexact uzawa
  algorithm for saddle point problems. SIAM Journal on Numerical Analysis
  34(3):1072--1092

\bibitem[{Cai et~al(2010)Cai, Cand{\`e}s, and Shen}]{cai2010singular}
Cai JF, Cand{\`e}s EJ, Shen Z (2010) A singular value thresholding algorithm
  for matrix completion. SIAM Journal on optimization 20(4):1956--1982

\bibitem[{Cai et~al(2018)Cai, Jiao, Zhuge, Tao, and Hou}]{CAI2018316}
Cai Y, Jiao Y, Zhuge W, et~al (2018) Partial multi-view spectral clustering.
  Neurocomputing 311:316--324

\bibitem[{Chao et~al(2022)Chao, Wang, Yang, Li, and Chu}]{chao2022incomplete}
Chao G, Wang S, Yang S, et~al (2022) Incomplete multi-view clustering with
  multiple imputation and ensemble clustering. Applied Intelligence pp 1--11

\bibitem[{Delaitre et~al(2010)Delaitre, Laptev, and
  Sivic}]{delaitre2010recognizing}
Delaitre V, Laptev I, Sivic J (2010) Recognizing human actions in still images:
  a study of bag-of-features and part-based representations. In: BMVC 2010-21st
  British Machine Vision Conference

\bibitem[{Deng et~al(2023)Deng, Wen, Liu, Yan, Xu, and Xu}]{deng2023projective}
Deng S, Wen J, Liu C, et~al (2023) Projective incomplete multi-view clustering.
  IEEE Transactions on Neural Networks and Learning Systems

\bibitem[{Fang et~al(2023)Fang, Li, Li, Gao, Jia, and
  Zhang}]{fang2023comprehensive}
Fang U, Li M, Li J, et~al (2023) A comprehensive survey on multi-view
  clustering. IEEE Transactions on Knowledge and Data Engineering

\bibitem[{Greene and Cunningham(2006)}]{greene2006practical}
Greene D, Cunningham P (2006) Practical solutions to the problem of diagonal
  dominance in kernel document clustering. In: Proceedings of the 23rd
  international conference on Machine learning, pp 377--384

\bibitem[{Greene and Cunningham(2013)}]{greene2013producing}
Greene D, Cunningham P (2013) Producing a unified graph representation from
  multiple social network views. In: Proceedings of the 5th annual ACM web
  science conference, pp 118--121

\bibitem[{Guo et~al(2023)Guo, Wang, Chi, Xu, Li, and Wu}]{guo2023scalable}
Guo W, Wang Z, Chi Z, et~al (2023) Scalable one-stage multi-view subspace
  clustering with dictionary learning. Knowledge-Based Systems 259:110092

\bibitem[{He et~al(2023)He, Zhang, and Wei}]{he2023scalable}
He Wj, Zhang Z, Wei Y (2023) Scalable incomplete multi-view clustering with
  adaptive data completion. Information Sciences p 119562

\bibitem[{Hu and Chen(2019)}]{hu2019doubly}
Hu M, Chen S (2019) Doubly aligned incomplete multi-view clustering. arXiv
  preprint arXiv:190302785

\bibitem[{Hu et~al(2021)Hu, Luo, Wang, Gao, Sun, and Yin}]{hu2021complete}
Hu Y, Luo C, Wang B, et~al (2021) Complete/incomplete multi-view subspace
  clustering via soft block-diagonal-induced regulariser. IET Computer Vision
  15(8):618--632

\bibitem[{Khan et~al(2023)Khan, Khan, Khan, Anwar, Ashraf, Atoum, Ahmad,
  Shahid, Ishrat, and Alghamdi}]{khan2023adaptive}
Khan MA, Khan GA, Khan J, et~al (2023) Adaptive weighted low-rank sparse
  representation for multi-view clustering. IEEE Access

\bibitem[{Li et~al(2022)Li, Tang, Zheng, Liu, Zhang, and Zhu}]{li2022high}
Li Z, Tang C, Zheng X, et~al (2022) High-order correlation preserved incomplete
  multi-view subspace clustering. IEEE Transactions on Image Processing
  31:2067--2080

\bibitem[{Lian et~al(2021)Lian, Xu, Wang, Li, Zhu, and Liu}]{lian2021partial}
Lian H, Xu H, Wang S, et~al (2021) Partial multiview clustering with locality
  graph regularization. International Journal of Intelligent Systems
  36(6):2991--3010

\bibitem[{Liang et~al(2022)Liang, Yang, and Xie}]{liang2022incomplete}
Liang N, Yang Z, Xie S (2022) Incomplete multi-view clustering with
  sample-level auto-weighted graph fusion. IEEE Transactions on Knowledge and
  Data Engineering 35(6):6504--6511

\bibitem[{Liu et~al(2012)Liu, Lin, Yan, Sun, Yu, and Ma}]{liu2012robust}
Liu G, Lin Z, Yan S, et~al (2012) Robust recovery of subspace structures by
  low-rank representation. IEEE transactions on pattern analysis and machine
  intelligence 35(1):171--184

\bibitem[{Liu et~al(2021{\natexlab{a}})Liu, Liu, Yang, Wang, and
  Zhou}]{liu2021hierarchical}
Liu J, Liu X, Yang Y, et~al (2021{\natexlab{a}}) Hierarchical multiple kernel
  clustering. In: Thirty-Fifth AAAI Conference on Artificial Intelligence,
  AAAI, pp 2--9

\bibitem[{Liu et~al(2021{\natexlab{b}})Liu, Liu, Zhang, Zhang, Tu, Wang, Zhou,
  Liang, Wang, and Yang}]{liu2021self}
Liu J, Liu X, Zhang Y, et~al (2021{\natexlab{b}}) Self-representation subspace
  clustering for incomplete multi-view data. In: Proceedings of the 29th ACM
  International Conference on Multimedia, pp 2726--2734

\bibitem[{Liu et~al(2023)Liu, Yang, Li, Li, and Xie}]{liu2023auto}
Liu M, Yang Z, Li L, et~al (2023) Auto-weighted collective matrix factorization
  with graph dual regularization for multi-view clustering. Knowledge-Based
  Systems 260:110145

\bibitem[{Liu and Lin(2023)}]{liu2023adaptive}
Liu SS, Lin L (2023) Adaptive weighted multi-view clustering. In: Conference on
  Health, Inference, and Learning, PMLR, pp 19--36

\bibitem[{Lu et~al(2012)Lu, Min, Zhao, Zhu, Huang, and
  Yan}]{10.1007/978-3-642-33786-4_26}
Lu CY, Min H, Zhao ZQ, et~al (2012) Robust and efficient subspace segmentation
  via least squares regression. In: Computer Vision -- ECCV 2012. Springer
  Berlin Heidelberg, Berlin, Heidelberg, pp 347--360

\bibitem[{Merino(1992)}]{merino1992topics}
Merino DI (1992) Topics in matrix analysis. The Johns Hopkins University

\bibitem[{Ren et~al(2018)Ren, Xiao, Xu, Guo, Chen, Wang, and
  Fang}]{ren2018robust}
Ren P, Xiao Y, Xu P, et~al (2018) Robust auto-weighted multi-view clustering.
  In: IJCAI, pp 2644--2650

\bibitem[{Samaria and Harter(1994)}]{samaria1994parameterisation}
Samaria FS, Harter AC (1994) Parameterisation of a stochastic model for human
  face identification. In: Proceedings of 1994 IEEE workshop on applications of
  computer vision, IEEE, pp 138--142

\bibitem[{Shi et~al(2022)Shi, Nie, Wang, and Li}]{shi2022self}
Shi S, Nie F, Wang R, et~al (2022) Self-weighting multi-view spectral
  clustering based on nuclear norm. Pattern Recognition 124:108429

\bibitem[{Shu et~al(2022)Shu, Zhang, Gao, Yang, Wang, and Gao}]{shu2022self}
Shu X, Zhang X, Gao Q, et~al (2022) Self-weighted anchor graph learning for
  multi-view clustering. IEEE Transactions on Multimedia

\bibitem[{Sun et~al(2021)Sun, Zhang, Wang, Zhou, Tu, Liu, Zhu, and
  Wang}]{sun2021scalable}
Sun M, Zhang P, Wang S, et~al (2021) Scalable multi-view subspace clustering
  with unified anchors. In: Proceedings of the 29th ACM International
  Conference on Multimedia, pp 3528--3536

\bibitem[{Tang et~al(2022)Tang, Cao, Zhang, and Jiang}]{tang2022consistent}
Tang K, Cao L, Zhang N, et~al (2022) Consistent auto-weighted multi-view
  subspace clustering. Pattern Analysis and Applications 25(4):879--890

\bibitem[{Vidal(2011)}]{vidal2011subspace}
Vidal R (2011) Subspace clustering. IEEE Signal Processing Magazine
  28(2):52--68

\bibitem[{Wan et~al(2023)Wan, Liu, Liu, Wang, Wen, Liang, Zhu, Liu, and
  Zhou}]{wan2023auto}
Wan X, Liu X, Liu J, et~al (2023) Auto-weighted multi-view clustering for
  large-scale data. arXiv preprint arXiv:230301983

\bibitem[{Wang et~al(2022)Wang, Wang, and Le}]{9914468}
Wang S, Wang Y, Le W (2022) Adaptive weight structure representation for
  multi-view subspace clustering. In: 2022 9th International Conference on
  Dependable Systems and Their Applications (DSA), pp 918--925,
  \doi{10.1109/DSA56465.2022.00129}

\bibitem[{Wen et~al(2018)Wen, Xu, and Liu}]{wen2018incomplete}
Wen J, Xu Y, Liu H (2018) Incomplete multiview spectral clustering with
  adaptive graph learning. IEEE transactions on cybernetics 50(4):1418--1429

\bibitem[{Wen et~al(2019)Wen, Zhang, Xu, Zhang, Fei, and Liu}]{wen2019unified}
Wen J, Zhang Z, Xu Y, et~al (2019) Unified embedding alignment with missing
  views inferring for incomplete multi-view clustering. In: Proceedings of the
  AAAI conference on artificial intelligence, pp 5393--5400

\bibitem[{Yin and Jiang(2023)}]{yin2023incomplete}
Yin J, Jiang J (2023) Incomplete multi-view clustering based on
  self-representation. Neural Processing Letters pp 1--15

\bibitem[{Zhang et~al(2015)Zhang, Fu, Liu, Liu, and Cao}]{zhang2015low}
Zhang C, Fu H, Liu S, et~al (2015) Low-rank tensor constrained multiview
  subspace clustering. In: Proceedings of the IEEE international conference on
  computer vision, pp 1582--1590

\bibitem[{Zhang et~al(2017)Zhang, Hu, Fu, Zhu, and Cao}]{zhang2017latent}
Zhang C, Hu Q, Fu H, et~al (2017) Latent multi-view subspace clustering. In:
  Proceedings of the IEEE conference on computer vision and pattern
  recognition, pp 4279--4287

\bibitem[{Zhang et~al(2020)Zhang, Wang, Hu, Cheng, Guo, Zhu, and
  Cai}]{zhang2020adaptive}
Zhang P, Wang S, Hu J, et~al (2020) Adaptive weighted graph fusion incomplete
  multi-view subspace clustering. Sensors 20(20):5755

\bibitem[{Zhao et~al(2016)Zhao, Liu, and Fu}]{zhao2016incomplete}
Zhao H, Liu H, Fu Y (2016) Incomplete multi-modal visual data grouping. In:
  IJCAI, pp 2392--2398

\bibitem[{Zhao et~al(2022{\natexlab{a}})Zhao, Zhang, Yang, and
  Chen}]{zhao2022incomplete}
Zhao L, Zhang J, Yang T, et~al (2022{\natexlab{a}}) Incomplete multi-view
  clustering based on weighted sparse and low rank representation. Applied
  Intelligence 52(13):14822--14838

\bibitem[{Zhao et~al(2023)Zhao, Yang, and Nie}]{zhao2023auto}
Zhao M, Yang W, Nie F (2023) Auto-weighted orthogonal and nonnegative graph
  reconstruction for multi-view clustering. Information Sciences 632:324--339

\bibitem[{Zhao et~al(2022{\natexlab{b}})Zhao, Dai, Wu, Peng, Liu, Bai, Tan,
  Wang, and Philip}]{zhao2022multi}
Zhao X, Dai Q, Wu J, et~al (2022{\natexlab{b}}) Multi-view tensor graph neural
  networks through reinforced aggregation. IEEE Transactions on Knowledge and
  Data Engineering 35(4):4077--4091

\bibitem[{Zhuge et~al(2017)Zhuge, Hou, Jiao, Yue, Tao, and
  Yi}]{zhuge2017robust}
Zhuge W, Hou C, Jiao Y, et~al (2017) Robust auto-weighted multi-view subspace
  clustering with common subspace representation matrix. PloS one
  12(5):e0176769

\end{thebibliography}



\end{document}